\documentclass[10pt]{article}
\usepackage{amstext}
\usepackage[fleqn]{amsmath}
\usepackage{amsfonts}
\usepackage{amssymb}
\usepackage{amsthm}
\usepackage{newlfont}
\usepackage{graphicx}
\usepackage{sectsty}
\theoremstyle{plain} \theoremstyle{theorem}
\newtheorem{theorem}{Theorem}[section]
\theoremstyle{example}

\theoremstyle{corollary}

\theoremstyle{lemma}

\theoremstyle{proposition}

\theoremstyle{axiom}

\theoremstyle{notation}

\theoremstyle{fact}

\theoremstyle{definition}
\newtheorem{definition}{Definition}[section]
\theoremstyle{remark}
\newtheorem{remark}{Remark}[section]
\numberwithin{equation}{section}
\begin{document}
%---------------------------------------------------------------------------------------------------------------------------------------------------------------------------------------
\title{An extended version of the $\;_{r+1}R_{s,k}(B,C,z)$ matrix function}
%---------------------------------------------------------------------------------------------------------------------------------------------------------------------------------------
%---------------------------------------------------------------------------------------------------------------------------------------------------------------------------------------
\author{Ayman Shehata \thanks{%
E-mail: drshehata2006@yahoo.com, aymanshehata@science.aun.edu.eg}\\
{\small Department of Mathematics, Faculty of Science, Assiut University, Assiut 71516, Egypt.}
}
%---------------------------------------------------------------------------------------------------------------------------------------------------------------------------------------
%---------------------------------------------------------------------------------------------------------------------------------------------------------------------------------------
\date{}
\maketitle{}
%---------------------------------------------------------------------------------------------------------------------------------------------------------------------------------------
%-------------------------------------------------------------------------------------------------------------------------------------------------------------------------------------
\begin{abstract}
%-------------------------------------------------------------------------------------------------------------------------------------------------------------------------------------
Recently, Shehata et al. \cite{skc} introduced the $\;_{r+1}R_{s}(B,C,z)$ matrix function and established some properties. The aim of this study established to devote and derive certain basic properties including analytic properties, recurrence matrix relations, differential properties, new integral representations, $k$-Beta transform, Laplace transform, fractional $k$-Fourier transform, fractional integral properties, the $k$-Riemann–Liouville and $k$-Weyl fractional integral and derivative operators an extended version of $\;_{r+1}R_{s,k}$ matrix function. We establish its relationships with other well known special matrix functions which have some particular cases in the context of three parametric Mittag-Leffler matrix function, $k$-Konhauser and $k$-Laguerre matrix polynomials. Finally, some special cases of the established formulas are also discussed.
%---------------------------------------------------------------------------------------------------------------------------------------------------------------------------------------
\end{abstract}
%-------------------------------------------------------------------------------------------------------------------------------------------------------------------------------------
%-------------------------------------------------------------------------------------------------------------------------------------------------------------------------------------

%---------------------------------------------------------------------------------------------------------------------------------------------------------------------------------------
\section{Introduction}
%-------------------------------------------------------------------------------------------------------------------------------------------------------------------------------------
Fractional calculus is the study of applications of derivatives and integrals of non-integer order. It is a generalized form of calculus, so it retains many properties of calculus. It is worth mentioning that, in recent times, theory of fractional calculus has developed quickly and played many important roles in science and engineering, serving as a powerful and very effective tool for various mathematical problems. It has been widely investigated in the last two decades. The hypergeometric function has a long history of mathematical and physical applications. They introduced integral representations of some $k$-confluent hypergeometric and $k$-hypergeometric functions. With the help of this new generalised pochhammer symbol. Joshi and Mittal \cite{jm}, Mubeen and Habibullah \cite{mh1, mh2},  Rahman et al. \cite{ram}, Sharma and Jain \cite{sj}, Zhou et al. \cite{zlx}, introduced an integral representation of $k$-Gamma and $k$Beta functions and some generalized $k$-hypergeometric functions. Mubeen et al. \cite{mr} also introduced $k$-analogue of Kummer's first formula and solution of some integral equations involving confluent $k$-hypergeometric functions. These studies were extended by Ahmad et al. \cite{aka}, Ali et al. \cite{aiih}, Diaz and Pariguan \cite{dp}, Farid et al. \cite{flaion}, Gupta and Bhatt \cite{gb}, Mittal and Joshi \cite{mj}, Mittal et al. \cite{mjs}, Romero and Cerutti \cite{rc}. Jain et al. \cite{jgoam}, Mubeen et al. \cite{mar}, Rahman et al. \cite{rg} introduced and derived some identities of $k$-Gamma matrix function, $k$-Beta matrix function, $k$-hypergeometric matrix functions and $k$-fractional integrations.

For the past four decades or so, in both mathematics and science, certain special matrix functions are crucial. J\'{o}dar and Sastre \cite{j1},  J\'{o}dar and Cort\'{e}s\cite{j2, j3, j4} researched the matrix analogues of the gamma, beta, and Gauss hypergeometric functions, which provided the basis for the special matrix functions  Bakhet et al.\cite{bjh}, Çekim et al. \cite{cdss}, Gezer and Kaanoglu \cite{gk}.The extended work of $\;_{r+1}R_{s}(P,Q,z)$ matrix functions is examined in some detail in  Sanjhira and Dave \cite{sd1}, Sanjhira and Dwivedi \cite{sd2}, Sanjhira et al. \cite{snd}, Shehata \cite{sa1, sa4, sa5}, Shehata et al. \cite{skc}, Varma et al. \cite{vct} for examples of several polynomials that have been introduced and investigated from a matrix perspective. The generalization of the $\;_{r+1}R_{s}(P,Q,z)$ function presented here was motivated by our investigations \cite{sd2, skc} of the properties of a class of polynomials which characterize is itself an interesting subject, the subsequent generalization of the $\;_{r+1}R_{s,k}(P,Q,z)$ matrix function would appear to be of mathematical interest on its own.
%---------------------------------------------------------------------------------------------------------------------------------------------------------------------------------------
\subsection{Preliminaries and some definitions}
%-------------------------------------------------------------------------------------------------------------------------------------------------------------------------------------
For this purpose, we will introduce the notations, properties, and definitions which we need in further sections. Throughout this paper, for a matrix $\Bbb{A}$ in $\Bbb{C}^{N\times N}$, its spectrum $\sigma(\Bbb{A})$ denotes the set of all eigenvalues of $\Bbb{A}$. The two-norm of an $\Bbb{C}^{N\times N}$ matrix $\Bbb{A}$ will be denoted by $||\Bbb{A}||_{2}$ and it is defined by (see \cite{j1, j2, j3, j4})
%-----------------------------------------------------------------------------------------------------------------------------------------------------------------
\begin{equation*}
\begin{split}
||\Bbb{A}||_{2}=\sup_{\zeta\neq 0}\frac{||\Bbb{A}\zeta||_{2}}{||\zeta||_{2}},
\end{split}
\end{equation*}
%-----------------------------------------------------------------------------------------------------------------------------------------------------------------
where for a vector $\zeta$ in $\Bbb{C}^{N}$, $||\zeta||_2=(\zeta^T\zeta)^\frac{1}{2}$ is the Euclidean norm of $\zeta$, then the $\zeta^T$ vector is the Hermitian transpose of
$\zeta$ (or, equivalently, the Hermitian transpose of the vector $\zeta$  that is viewed as a matrix) .
%-----------------------------------------------------------------------------------------------------------------------------------------------------------------

Let us denote the real numbers $\Bbb{M}(\Bbb{A})$ and $\mathbf{m}(\Bbb{A})$ as in the following
%-----------------------------------------------------------------------------------------------------------------------------------------------------------------
\begin{equation}
\begin{split}
\Bbb{M}(\Bbb{A})=\max\{Re(\zeta): \zeta\in \sigma(\Bbb{A})\};\quad \mathbf{m}(\Bbb{A})=\min\{Re(\zeta): \zeta\in \sigma(\Bbb{A})\}.\label{1.1}
\end{split}
\end{equation}
%-----------------------------------------------------------------------------------------------------------------------------------------------------------------
If $\mathbf{\Upsilon}(\zeta)$ and $\mathbf{\Psi}(\zeta)$ are holomorphic functions of the complex variable $\zeta$, which are defined in an open set $\Omega$ of the complex plane, and $\Bbb{A}$, $\Bbb{B}$ are matrices in $\Bbb{C}^{N\times N}$ with $\sigma(\Bbb{A})\subset\Omega$ and $\sigma(\Bbb{B}) \subset \Omega$, such that $\Bbb{A}\Bbb{B}=\Bbb{B}\Bbb{A}$, then the properties of the matrix functional calculus in \cite{ds}, it follows that
%-----------------------------------------------------------------------------------------------------------------------------------------------------------------
\begin{equation*}
\begin{split}
\mathbf{\Upsilon}(\Bbb{A})\mathbf{\Psi}(\Bbb{B})=\mathbf{\Psi}(\Bbb{B})\mathbf{\Upsilon}(\Bbb{A}).
\end{split}
\end{equation*}
%-----------------------------------------------------------------------------------------------------------------------------------------------------------------
%-------------------------------------------------------------------------------------------------------------------------------------------------------------------------------------
\begin{definition}  \cite{mar}
%-------------------------------------------------------------------------------------------------------------------------------------------------------------------------------------
For $k>0$, $\Bbb{A}$ and $\Bbb{B}$ are positive stable matrices in $\Bbb{C}^{N\times N}$, then the $\kappa$-Gamma and $\kappa$-Beta matrix functions are defined by (see  \cite{kac, mra})
%-----------------------------------------------------------------------------------------------------------------------------------------------------------------
\begin{equation}
\begin{split}
\Gamma_{\kappa}(\Bbb{A})=\int_{0}^{\infty}t^{\Bbb{A}-I}e^{-\frac{t^{\kappa}}{\kappa}}dt=\kappa^{\frac{\Bbb{A}}{\kappa}-1}\Gamma(\frac{\Bbb{A}}{\kappa})\label{1.2}
\end{split}
\end{equation}
%-----------------------------------------------------------------------------------------------------------------------------------------------------------------
and
%-----------------------------------------------------------------------------------------------------------------------------------------------------------------
\begin{equation}
\begin{split}
\Bbb{B}_{\kappa}(\Bbb{A},\Bbb{B})=\frac{1}{\kappa}\int_{0}^{1}t^{\frac{\Bbb{A}}{\kappa}-I}(1-t)^{\frac{\Bbb{B}}{\kappa}-I}dt=\Gamma_{\kappa}(\Bbb{A})\Gamma_{\kappa}(\Bbb{B})\Gamma_{\kappa}^{-1}(\Bbb{A}+\Bbb{B}),\label{1.3}
\end{split}
\end{equation}
%-----------------------------------------------------------------------------------------------------------------------------------------------------------------
where $I$ is the identity matrix in $\Bbb{C}^{N\times N}$.
%-------------------------------------------------------------------------------------------------------------------------------------------------------------------------------------
\end{definition}
%-------------------------------------------------------------------------------------------------------------------------------------------------------------------------------------
%-----------------------------------------------------------------------------------------------------------------------------------------------------------------
Furthermore, if $\Bbb{A}$ is a matrix such that (see \cite{mar, rg})
%-----------------------------------------------------------------------------------------------------------------------------------------------------------------
\begin{equation}
\begin{split}
\Bbb{A}+k\ell I\quad \text{is\; an invertible matrix\; for\; all\; integers}\; \ell\geq0, \forall k>0\label{1.4}
\end{split}
\end{equation}
%-----------------------------------------------------------------------------------------------------------------------------------------------------------------
then the $k$-pochammer matrix symbol is defined as
%-------------------------------------------------------------------------------------------------------------------------------------------------------------------------------------
\begin{equation}
\begin{split}
(\Bbb{A})_{n,\kappa}=\Bbb{A}(\Bbb{A}+\kappa I)\ldots(\Bbb{A}+(n-1)\kappa I)=\Gamma_{\kappa}{(\Bbb{A}+n\kappa I)}\Gamma^{-1}_{\kappa}{(\Bbb{A})}\;;
n\geq 1\;, (\Bbb{A})_{0,k}=I; (\Bbb{A})_{n,1}=(\Bbb{A})_{n}.\label{1.5}
\end{split}
\end{equation}
%-------------------------------------------------------------------------------------------------------------------------------------------------------------------------------------
Also, they have provided some useful results
%-------------------------------------------------------------------------------------------------------------------------------------------------------------------------------------
\begin{equation}
\begin{split}
(\Bbb{A})_{n\ell,\kappa}=\ell^{n\ell}\bigg{(}\frac{\Bbb{A}}{\ell}\bigg{)}_{n,\kappa}\bigg{(}\frac{\Bbb{A}+kI}{\ell}\bigg{)}_{n,\kappa}\ldots\bigg{(}\frac{\Bbb{A}+k(\ell-1)I}{\ell}\bigg{)}_{n,\kappa}\label{1.6}
\end{split}
\end{equation}
%-------------------------------------------------------------------------------------------------------------------------------------------------------------------------------------
and
%-------------------------------------------------------------------------------------------------------------------------------------------------------------------------------------
\begin{equation}
\begin{split}
(1-k\zeta)^{-\frac{\Bbb{A}}{k}}=\sum_{n=0}^{\infty}\frac{(\Bbb{A})_{n,\kappa}}{n!}\zeta^{n}.\label{1.7}
\end{split}
\end{equation}
%-------------------------------------------------------------------------------------------------------------------------------------------------------------------------------------
If $\Phi(\kappa,\ell)$ is matrix in $\Bbb{C}^{N\times N}$ for $\kappa\geq 0$, $\ell\geq 0$, then in an analogous way to the proof of Lemma
\textbf{11} \cite{j2}, it follows that
%-------------------------------------------------------------------------------------------------------------------------------------------------------------------------------------
\begin{equation}
\begin{split}
\sum_{\kappa=0}^{\infty}\sum_{\ell=0}^{\infty}\Phi(\kappa,\ell)=\sum_{\kappa=0}^{\infty}\sum_{\ell=0}^{[\frac{1}{2}\kappa]}\Phi(\kappa-2\ell,\ell),\\
\sum_{\kappa=0}^{\infty}\sum_{\ell=0}^{\infty}\Phi(\kappa,\ell)=\sum_{\kappa=0}^{\infty}\sum_{\ell=0}^{\kappa}\Phi(\kappa-\ell,\ell).\label{1.8}
\end{split}
\end{equation}
%-------------------------------------------------------------------------------------------------------------------------------------------------------------------------------------
Similarly to (\ref{1.8}), we can write
%-----------------------------------------------------------------------------------------------------------------------------------------------------------------
\begin{equation}
\begin{split}
\sum_{\kappa=0}^{\infty}\sum_{\ell=0}^{[\frac{1}{2}\kappa]}\Phi(\kappa,\ell)=\sum_{\kappa=0}^{\infty}\sum_{\ell=0}^{\infty}\Phi(\kappa+2\ell,\ell),\\
\sum_{\kappa=0}^{\infty}\sum_{\ell=0}^{\kappa}\Phi(\kappa,\ell)=\sum_{\kappa=0}^{\infty}\sum_{\ell=0}^{\infty}\Phi(\kappa+\ell,\ell).\label{1.9}
\end{split}
\end{equation}
%-------------------------------------------------------------------------------------------------------------------------------------------------------------------------------------
\begin{definition}
%-------------------------------------------------------------------------------------------------------------------------------------------------------------------------------------
%-----------------------------------------------------------------------------------------------------------------------------------------------------------------
The hypergeometric matrix function $_{2}F_{1,\kappa}(\Bbb{A},\Bbb{B};C;\zeta)$ has been given in the form (\cite{mar, rg})
%-----------------------------------------------------------------------------------------------------------------------------------------------------------------
\begin{equation}
\begin{split}
_{2}F_{1,\kappa}(\Bbb{A},\Bbb{B};\Bbb{C};\zeta)=_{2}F_{1}((\Bbb{A},\kappa),(\Bbb{B},\kappa);(C,\kappa);\zeta)=\sum_{\ell=0}^{\infty}\frac{(\Bbb{A})_{\ell,\kappa}(\Bbb{B})_{\ell,\kappa}[(C)_{\ell,\kappa}]^{-1}}{n!}\zeta^\kappa,\label{1.10}
\end{split}
\end{equation}
%-----------------------------------------------------------------------------------------------------------------------------------------------------------------
for matrices $\Bbb{A}$, $\Bbb{B}$ and $C$ in $\Bbb{C}^{N\times N}$ such that $C+\ell I$ is an invertible matrix for all integers $\ell\geq 0$ and for $|\zeta|<1$. It has been
seen by J\'{o}dar and Cort\'{e}s \cite{j2} that the series is absolutely convergent for $|\zeta|=1$ when
%-----------------------------------------------------------------------------------------------------------------------------------------------------------------
\begin{equation*}
\begin{split}
\mathbf{m}(\Bbb{C})> \Bbb{M}(\Bbb{A})+\Bbb{M}(\Bbb{B}),
\end{split}
\end{equation*}
%-----------------------------------------------------------------------------------------------------------------------------------------------------------------
where $\mathbf{m}(\Bbb{A})$ and $\Bbb{M}(\Bbb{A})$ in (\ref{1.1}) for any matrix $\Bbb{A}$ in $\Bbb{C}^{N\times N}$.
%-------------------------------------------------------------------------------------------------------------------------------------------------------------------------------------
\end{definition}
%-------------------------------------------------------------------------------------------------------------------------------------------------------------------------------------
%-------------------------------------------------------------------------------------------------------------------------------------------------------------------------------------
\begin{definition}
%-------------------------------------------------------------------------------------------------------------------------------------------------------------------------------------
For $m$ and $n$ are finite positive integers, the $\;_{m}R_{n}$ matrix function defined as (see \cite{sd2, skc})
%-----------------------------------------------------------------------------------------------------------------------------------------------------------------
\begin{equation}
\begin{split}
&\;_{m}R_{n}(A_{1},A_{2},\ldots,A_{m};Q_{1},Q_{2},\ldots,Q_{n};P,Q;\zeta)\\
=&\sum_{k=0}^{\infty}\frac{\zeta^{k}}{k!}(A_{1})_{k}(A_{2})_{k}\ldots(A_{m})_{k}[(Q_{1})_{k}]^{-1}[(Q_{2})_{k}]^{-1}\ldots[(Q_{n})_{k}]^{-1}\Gamma^{-1}(\ell P+Q)\\
=&\sum_{k=0}^{\infty}\frac{\zeta^{k}}{k!}\prod_{i=1}^{m}(A_{i})_{k}\bigg{[}\prod_{j=1}^{n}(Q_{i})_{k}\bigg{]}^{-1}\Gamma^{-1}(\ell P+Q),\label{1.11}
\end{split}
\end{equation}
%-----------------------------------------------------------------------------------------------------------------------------------------------------------------
where $A_{i}$ ; $1\leq i\leq m$ and $Q_{j}$; $1\leq j\leq n$ are matrices in $\Bbb{C}^{N\times N}$ such that
%-----------------------------------------------------------------------------------------------------------------------------------------------------------------
\begin{equation}
\begin{split}
Q_{j}+kI \quad \text{are\; invertible\; matrices\; for\; all\; integers}\; k\geq0.\label{1.12}
\end{split}
\end{equation}
%-----------------------------------------------------------------------------------------------------------------------------------------------------------------
%-------------------------------------------------------------------------------------------------------------------------------------------------------------------------------------
\end{definition}
%-------------------------------------------------------------------------------------------------------------------------------------------------------------------------------------
%-------------------------------------------------------------------------------------------------------------------------------------------------------------------------------------
\begin{definition}
%-------------------------------------------------------------------------------------------------------------------------------------------------------------------------------------
 The $k$-Riemann Liouville fractional $k$-integral and derivative operators of order $\mu$ defined as follows (Mubeen and Habibullah see \cite{mh1, mh2})
%-----------------------------------------------------------------------------------------------------------------------------------------------------------------
\begin{equation}
\begin{split}
I_{a,k}^{\mu}\mathbf{\Psi}(x)=\frac{1}{k\Gamma_{k}(\mu)}\int_{a}^{x}(x-t)^{\frac{\mu}{k}-1}\mathbf{\Psi}(t)dt,\mu\in R^{+}\label{1.13}
\end{split}
\end{equation}
%-----------------------------------------------------------------------------------------------------------------------------------------------------------------
and
%-----------------------------------------------------------------------------------------------------------------------------------------------------------------
%---------------------------------------------------------------------------------------------------------------------------------------------------------------------------------------------
\begin{equation}
\begin{split}
\bigg{(}\mathbb{D}_{a^{+},k}^{\mu}\mathbf{\Psi}\bigg{)}(x)=\bigg{(}\frac{d}{dx}\bigg{)}^{n}\bigg{(}\mathbb{I}_{a^{+},k}^{n-\mu}\mathbf{\Psi}\bigg{)}(x), \label{1.14}
\end{split}
\end{equation}
%---------------------------------------------------------------------------------------------------------------------------------------------------------------------------------------------
where $Re(\mu)>0$ and $k$ be any positive real number.
%-------------------------------------------------------------------------------------------------------------------------------------------------------------------------------------
\end{definition}
%-------------------------------------------------------------------------------------------------------------------------------------------------------------------------------------
\begin{definition}
%-------------------------------------------------------------------------------------------------------------------------------------------------------------------------------------
For $\alpha\in C$, $Re(\alpha)>0$ and $k\in R^{+}$, $f$ belonging to $S(R)$, the Weyl fractional $k$-integral operator and $k$-Weyl fractional derivative are defined as (see \cite{rc, rcd})
%-----------------------------------------------------------------------------------------------------------------------------------------------------------------
\begin{equation}
\begin{split}
\bigg{[}\mathbb{W}_{k}^{\alpha}\mathbf{\Psi}\bigg{]}(x)=\frac{1}{k\Gamma_{k}(\alpha)}\int_{x}^{\infty}(t-x)^{\frac{\alpha}{k}-1}\mathbf{\Psi}(t)dt\label{1.15}
\end{split}
\end{equation}
%---------------------------------------------------------------------------------------------------------------------------------------------------------------------------------------------
and
%---------------------------------------------------------------------------------------------------------------------------------------------------------------------------------------------
\begin{equation}
\begin{split}
\bigg{[}\mathbb{W}_{k}^{-\alpha}\mathbf{\Psi}\bigg{]}(x)=-\frac{d}{dx}\bigg{[}\mathbb{W}_{k}^{1-\alpha}\mathbf{\Psi}\bigg{]}(x).\label{1.16}
\end{split}
\end{equation}
%---------------------------------------------------------------------------------------------------------------------------------------------------------------------------------------------
%-------------------------------------------------------------------------------------------------------------------------------------------------------------------------------------
\end{definition}
%-------------------------------------------------------------------------------------------------------------------------------------------------------------------------------------
%-------------------------------------------------------------------------------------------------------------------------------------------------------------------------------------
\begin{definition}
%-------------------------------------------------------------------------------------------------------------------------------------------------------------------------------------
%---------------------------------------------------------------------------------------------------------------------------------------------------------------------------------------------
The $\mathfrak{B}[f(t);A,B]$ $k$-Beta transform, $\mathfrak{L}[f(t);\mathbf{s}]$ Laplace transform and $\mathfrak{F}[f(t)]$ Fractional $k$-Fourier transform of $f(t)$ are defined by
%---------------------------------------------------------------------------------------------------------------------------------------------------------------------------------------------
\begin{equation}
\begin{split}
\mathfrak{B}[\mathbf{\Psi}(t);A,B]=\frac{1}{k}\int_{0}^{1}t^{\frac{A}{k}-I}(1-t)^{\frac{B}{k}-I}\mathbf{\Psi}(t)dt.\label{1.17}
\end{split}
\end{equation}
%---------------------------------------------------------------------------------------------------------------------------------------------------------------------------------------------
%---------------------------------------------------------------------------------------------------------------------------------------------------------------------------------------------
\begin{equation}
\begin{split}
\mathfrak{L}[\mathbf{\Psi}(t);\mathbf{s}]=\int_{0}^{\infty}e^{-\mathbf{s}t}\mathbf{\Psi}(t)dt=F(\mathbf{s}),\mathbf{s}\in\mathcal{C}.\label{1.18}
\end{split}
\end{equation}
%---------------------------------------------------------------------------------------------------------------------------------------------------------------------------------------------
and
%---------------------------------------------------------------------------------------------------------------------------------------------------------------------------------------------
\begin{equation}
\begin{split}
\mathfrak{F}[\mathbf{\Psi}(t)]=\int_{-\infty}^{0}e^{iw^{\frac{1}{\alpha}}z}\mathbf{\Psi}(z)dz.\label{1.19}
\end{split}
\end{equation}
%-------------------------------------------------------------------------------------------------------------------------------------------------------------------------------------
\end{definition}
%-------------------------------------------------------------------------------------------------------------------------------------------------------------------------------------
%-------------------------------------------------------------------------------------------------------------------------------------------------------------------------------------
The motive of the current study is to investigate the analytical and fractional integral and derivative properties of $\;_{r}R_{s,k}$ matrix function, as well as to emphasize the importance of their applications in diverse research areas. This function is an amalgamation of generalized Mittag–Leffler function and generalized hypergeometric function which plays an important role in the theory of mathematical analysis, fractional calculus and statistics and has significant applications in the field of free electron laser equations and fractional kinetic equations. For literature survey of fractional integral operators, researchers can refer to the papers of Farid et al. \cite{flaion}, Jain et al. \cite{jgoam}, Kilbas et al. \cite{kst}. In the present sequel to the aforementioned and many other recent investigations. In Section 2, some properties, recurrence matrix relations, differential properties, new integral representations, $k$-Beta transform, Laplace Transform, Fractional $k$-Fourier transform, fractional integral properties, the $k$-Riemann–Liouville and $k$-Weyl fractional integral and derivative operators of the extended $\;_{r}R_{s,k}$ matrix functions are established and discussed. In section 3, we derive some properties of $\;_{r}R_{s}$ matrix functions. Finally, we give a definition of the novel generalized type $\;_{r}R_{s,k}$ matrix functions.
%-----------------------------------------------------------------------------------------------------------------------------------------------------------------
\section{Main Results}
%-----------------------------------------------------------------------------------------------------------------------------------------------------------------
For our aim section, we establish and discuss the properties convergence, recurrence matrix relations, $k$-integrals representations and differential properties, $k$-fractional operators, $k$-Beta, Laplace, Fractional $k$-Fourier transforms, fractional integral and derivative operators for the $\;_{r+1}R_{s,k}$ matrix function.
%---------------------------------------------------------------------------------------------------------------------------------------------------------------------------------------------
%-----------------------------------------------------------------------------------------------------------------------------------------------------------------
\begin{definition}
%-----------------------------------------------------------------------------------------------------------------------------------------------------------------
Let $A$, $B$, $C$, $P_{i}$, $1\leq i\leq r$ and $Q_{j}$; $1\leq j\leq s$ be matrices in $\Bbb{C}^{N\times N}$, $Re(B)>0$, $Re(C)>0$, $Re(P_{i})>0$, $Re(Q_{j})>0$, $k\in R^{+}$,
then we define the $\;_{r+1}R_{s,k}(B,C,z)$ matrix function as
%-----------------------------------------------------------------------------------------------------------------------------------------------------------------
\begin{equation}
\begin{split}
&\;_{r+1}R_{s,k}(A,P_{1},P_{2},\ldots,P_{r};Q_{1},Q_{2},\ldots,Q_{s};B,C;z)\\
=&\sum_{\ell=0}^{\infty}\frac{z^{\ell}}{\ell !}(A)_{\ell,k}\prod_{i=1}^{r}(P_{i})_{\ell,k}\bigg{[}\prod_{j=1}^{s}(Q_{j})_{\ell,k}\bigg{]}^{-1}\Gamma^{-1}_{k}(\ell B+C), \label{2.1}
\end{split}
\end{equation}
%-----------------------------------------------------------------------------------------------------------------------------------------------------------------
for $r$ and $s$ are finite positive integers, such that
%-----------------------------------------------------------------------------------------------------------------------------------------------------------------
\begin{equation}
\begin{split}
Q_{j}+k\ell I \quad \text{are\; invertible\; matrices\; for\; all\; integers}\; \ell \geq0,\label{2.2}
\end{split}
\end{equation}
%---------------------------------------------------------------------------------------------------------------------------------------------------------------------------------------
\end{definition}
%-----------------------------------------------------------------------------------------------------------------------------------------------------------------
Now, we investigate the convergence of the following series, one gets
%-----------------------------------------------------------------------------------------------------------------------------------------------------------------
%----------------------------------------------------------------------------------------------------------------------------------------------------------------------------------------------
\begin{eqnarray*}
\begin{split}
\frac{1}{R}&=\limsup_{\ell\rightarrow\infty}\left(\|U_{\ell}\|\right)^{\frac{1}{\ell}}=\lim_{\ell\to\infty}\sup\left(\bigg{\|}\frac{(A)_{\ell,k}\prod_{i=1}^{r}(P_{i})_{\ell,k}[\prod_{j=1}^{s,k}(Q_{j})_{\ell,k}]^{-1}\Gamma_{k}^{-1}(\ell B+C)}{\ell !}\bigg{\|}\right)^{\frac{1}{\ell}}\\
&=\limsup_{\ell\rightarrow\infty}\bigg{\|}\frac{\Gamma_{k}(A+k\ell I)\Gamma_{k}^{-1}(A)\prod_{i=1}^{r}\Gamma_{k}(P_{i}+k\ell I)\Gamma_{k}^{-1}(P_{i})\prod_{j=1}^{s}\Gamma_{k}^{-1}(Q_{j}+k\ell I)\Gamma_{k}(Q_{j})\Gamma_{k}^{-1}(\ell B+C)}{\ell !}\bigg{\|}^{\frac{1}{\ell}}\\
&=\limsup_{\ell\rightarrow\infty}\bigg{\|}\sqrt{2\pi}e^{-(A+k\ell I)}(A+k\ell I)^{A+k\ell I-\frac{1}{2}I}\bigg{(}\sqrt{2\pi}e^{-A}(A)^{A-\frac{1}{2}I}\bigg{)}^{-1}\\
&\times\prod_{i=1}^{r}\sqrt{2\pi}e^{-(P_{i}+k\ell I)}(P_{i}+\ell I)^{P_{i}+k\ell I-\frac{1}{2}I}\bigg{(}\sqrt{2\pi}e^{-P_{i}}(P_{i})^{P_{i}-\frac{1}{2}I}\bigg{)}^{-1}\\
&\times\prod_{j=1}^{s}\frac{1}{\sqrt{2\pi}}e^{(Q_{j}+k\ell I)}(Q_{j}+k\ell I)^{-Q_{j}-k\ell I+\frac{1}{2}I}\bigg{(}\frac{1}{\sqrt{2\pi}}e^{(Q_{j})}(Q_{j})^{-Q_{j}+\frac{1}{2}I}\bigg{)}\\
&\times\frac{1}{\sqrt{2\pi}}e^{(\ell B+C)}(\ell B+C)^{-\ell B-C+\frac{1}{2}I}\frac{1}{\sqrt{2\pi}e^{-\ell-1}\ell^{\ell+\frac{1}{2}}}\bigg{\|}^{\frac{1}{\ell}}\\
&\thickapprox\limsup_{\ell\rightarrow\infty}\bigg{\|}e^{-(A+k\ell I)}(A+k\ell I)^{A+k\ell I-\frac{1}{2}I}\prod_{i=1}^{r}e^{-(P_{i}+k\ell I)}(P_{i}+k\ell I)^{P_{i}+k\ell I-\frac{1}{2}I}\\
&\times\prod_{j=1}^{s}e^{(Q_{j}+k\ell I)}(Q_{j}+k\ell I)^{-Q_{j}-k\ell I+\frac{1}{2}I}e^{(\ell B+C)}(\ell B+C)^{-\ell B-C+\frac{1}{2}I}\frac{1}{e^{-\ell-1}\ell^{\ell+\frac{1}{2}}}\bigg{\|}^{\frac{1}{\ell}}
\end{split}
\end{eqnarray*}
%---------------------------------------------------------------------------------------------------------------------------------------------------------------------------------------------
%----------------------------------------------------------------------------------------------------------------------------------------------------------------------------------------------
\begin{eqnarray*}
\begin{split}
&\thickapprox\limsup_{\ell\rightarrow\infty}\bigg{\|}(A+k\ell I)^{A+k\ell I-\frac{1}{2}I}\prod_{i=1}^{r}\prod_{j=1}^{s}e^{-(A+k\ell I)-(P_{i}+k\ell I)+Q_{j}+k\ell I+\ell B+C+\ell+1}(P_{i}+k\ell I)^{P_{i}+k\ell I-\frac{1}{2}I}\\
&\times(Q_{j}+k\ell I)^{-Q_{j}-k\ell I+\frac{1}{2}I}(\ell B+C)^{-\ell B-C+\frac{1}{2}I}\ell^{-\ell-\frac{1}{2}}\bigg{\|}^{\frac{1}{\ell}}\\
&\thickapprox\limsup_{\ell\rightarrow\infty}\bigg{\|}(A+k\ell I)^{A+k\ell I-\frac{1}{2}I}\prod_{i=1}^{r}\prod_{j=1}^{s}e^{Q_{j}+C+\ell B-k\ell I +\ell I+I-A-P_{i}}(P_{i}+k\ell I)^{P_{i}+k\ell I-\frac{1}{2}I}\\
&\times(Q_{j}+k\ell I)^{-Q_{j}-k\ell I+\frac{1}{2}I}(\ell B+C)^{-\ell B-C+\frac{1}{2}I}\ell^{-\ell-\frac{1}{2}}\bigg{\|}^{\frac{1}{\ell}}\\
&\thickapprox\limsup_{\ell\rightarrow\infty}\bigg{\|}(A+k\ell I)^{A+k\ell I-\frac{1}{2}I}\prod_{i=1}^{r}\prod_{j=1}^{s}e^{(B+(1-k)I)\ell}(P_{i}+k\ell I)^{P_{i}+k\ell I-\frac{1}{2}I}\\
&\times(Q_{j}+k\ell I)^{-Q_{j}-k\ell I+\frac{1}{2}I}(\ell B+C)^{-\ell B-C+\frac{1}{2}I}\ell^{-\ell-\frac{1}{2}}\bigg{\|}^{\frac{1}{\ell}}\\
&\thickapprox {\|}e^{B+(1-k)I}{\|}\limsup_{\ell\rightarrow\infty}\bigg{\|}\prod_{i=1}^{r}\prod_{j=1}^{s}\frac{(A+k\ell I)(P_{i}+k\ell I)(Q_{j}+k\ell I)^{-1}(\ell B+C)^{-B}}{\ell}\bigg{\|}^{k}\\
&\times\bigg{\|}(A+k\ell I)^{A-\frac{1}{2}I}(P_{i}+k\ell I)^{P_{i}-\frac{1}{2}I}(Q_{j}+k\ell I)^{-Q_{j}+\frac{1}{2}I}(\ell B+C)^{-C+\frac{1}{2}I}\ell^{-\frac{1}{2}}\bigg{\|}^{\frac{1}{\ell}}.
\end{split}
\end{eqnarray*}
%-----------------------------------------------------------------------------------------------------------------------------------------------------------------
The above limit shows that the following:
%---------------------------------------------------------------------------------------------------------------------------------------------------------------------------------------------
\begin{enumerate}
%---------------------------------------------------------------------------------------------------------------------------------------------------------------------------------------------
\item If $r>s+1$, then the series in (\ref{2.1}) diverges for $z\neq0$.
%---------------------------------------------------------------------------------------------------------------------------------------------------------------------------------------------
\item If $r\leq s$, then the series in (\ref{2.1}) converges for all finite $z$.
%---------------------------------------------------------------------------------------------------------------------------------------------------------------------------------------------
\item If $r=s+1$, then the series in (\ref{2.1}) converges for all $|z|<\frac{1}{k}$ and diverges for all $|z|>\frac{1}{k}$.
%---------------------------------------------------------------------------------------------------------------------------------------------------------------------------------------------
\item If $r=s+1$, then the series in (\ref{2.1}) is absolutely convergent on the circle $|z|=\frac{1}{k}$ when
%---------------------------------------------------------------------------------------------------------------------------------------------------------------------------------------------
\begin{equation*}
\begin{split}
\sum_{j=1}^{s}\mathbf{m}(Q_{j})>\sum_{i=1}^{r}\Bbb{M}(P_{i})+\Bbb{M}(A).
\end{split}
\end{equation*}
%---------------------------------------------------------------------------------------------------------------------------------------------------------------------------------------------
\item If $r=s+1$, then the series (\ref{2.1}) is diverges for $|z|=\frac{1}{k}$ when
%---------------------------------------------------------------------------------------------------------------------------------------------------------------------------------------------
\begin{equation*}
\begin{split}
\sum_{j=0}^{s}\mathbf{m}(Q_{j})\leq\sum_{i=0}^{r}\Bbb{M}(P_{i})+\Bbb{M}(A)-k,
\end{split}
\end{equation*}
%-----------------------------------------------------------------------------------------------------------------------------------------------------------------
\item If $r=s+1$, then the series (\ref{2.1}) is conditionally convergent for $|z|=\frac{1}{k}$ when
%---------------------------------------------------------------------------------------------------------------------------------------------------------------------------------------------
\begin{equation*}
\begin{split}
\sum_{i=0}^{r}\Bbb{M}(P_{i})+\Bbb{M}(A)-k<\sum_{j=0}^{s}m(Q_{j})\leq\sum_{i=0}^{r}\Bbb{M}(P_{i})+\Bbb{M}(A).
\end{split}
\end{equation*}
%---------------------------------------------------------------------------------------------------------------------------------------------------------------------------------------------
where $\Bbb{M}(P_{i})$; $1\leq i\leq r$  and $\mathbf{m}(Q_{j})$; $1\leq j\leq s$ are defined in (\ref{1.1}).
%---------------------------------------------------------------------------------------------------------------------------------------------------------------------------------------------
\end{enumerate}
%----------------------------------------------------------------------------------------------------------------------------------------------------------------------------------------------
%---------------------------------------------------------------------------------------------------------------------------------------------------------------------------------------------
\begin{theorem}
%-----------------------------------------------------------------------------------------------------------------------------------------------------------------
The following relations for $\;_{r+1}R_{s,k}$ hold true
%-----------------------------------------------------------------------------------------------------------------------------------------------------------------
\begin{equation}
\begin{split}
&(A-P_{i})\;_{r+1}R_{s,k}=A\;_{r+1}R_{s,k}(A+kI)-P_{i}\;_{r+1}R_{s,k}(P_{i}+k I),i=1,2,3,\ldots,r,\\
&(P_{\nu}-P_{i})\;_{r+1}R_{s,k}=P_{\nu}\;_{r+1}R_{s,k}(P_{\nu}+kI)-P_{i}\;_{r+1}R_{s,k}(P_{i}+k I)\,;i\neq \nu;i,\nu=1,2,3,\ldots,r,\label{2.3}
\end{split}
\end{equation}
%-----------------------------------------------------------------------------------------------------------------------------------------------------------------
%-----------------------------------------------------------------------------------------------------------------------------------------------------------------
\begin{equation}
\begin{split}
(Q_{\nu}-Q_{j})\;_{r+1}R_{s,k}=(Q_{\nu}-kI)\;_{r+1}R_{s,k}(Q_{\nu}-kI)-(Q_{j}-kI)\;_{r+1}R_{s,k}(Q_{j}-kI)\,;\nu\neq j;\nu,j=1,2,\ldots,s\label{2.4}
\end{split}
\end{equation}
%-----------------------------------------------------------------------------------------------------------------------------------------------------------------
and
%-----------------------------------------------------------------------------------------------------------------------------------------------------------------
\begin{equation}
\begin{split}
&(A-Q_{j}+kI)\;_{r+1}R_{s,k}=A\;_{r+1}R_{s,k}(A+kI)-(Q_{j}-kI)\;_{r+1}R_{s,k}(Q_{j}-kI)\,;j=1,2,\ldots,s,\\
&(P_{i}-Q_{j}+kI)\;_{r+1}R_{s,k}=P_{i}\;_{r+1}R_{s,k}(P_{i}+kI)-(Q_{j}-kI)\;_{r+1}R_{s,k}(Q_{j}-kI)\,;i=1,2,3,\ldots,r,\,j=1,2,\ldots,s.\label{2.5}
\end{split}
\end{equation}
%---------------------------------------------------------------------------------------------------------------------------------------------------------------------------------------------
\end{theorem}
%-----------------------------------------------------------------------------------------------------------------------------------------------------------------
%---------------------------------------------------------------------------------------------------------------------------------------------------------------------------------------------
\begin{proof}
%---------------------------------------------------------------------------------------------------------------------------------------------------------------------------------------------
By using the following property
%-----------------------------------------------------------------------------------------------------------------------------------------------------------------
\begin{equation*}
\begin{split}
A(A+kI)_{\ell,k}=(A+k\ell I)(A)_{\ell,k},
\end{split}
\end{equation*}
%---------------------------------------------------------------------------------------------------------------------------------------------------------------------------------------------
%-----------------------------------------------------------------------------------------------------------------------------------------------------------------
we get the matrix contiguous function relation
%-----------------------------------------------------------------------------------------------------------------------------------------------------------------
\begin{equation}
\begin{split}
&A\;_{r+1}R_{s,k}(A+kI)=A\sum_{\ell=0}^{\infty}\frac{z^{\ell}}{\ell\;!}(A+kI)_{\ell,k}\prod_{i=1}^{r}(P_{i})_{\ell,k}\bigg{[}\prod_{j=1}^{s}(Q_{j})_{\ell,k}\bigg{]}^{-1}\Gamma^{-1}_{k}(\ell B+C)\\
&=\sum_{\ell=0}^{\infty}\frac{z^{\ell}}{\ell !}(A+k\ell I)(A)_{\ell,k}\prod_{i=1}^{r}(P_{i})_{\ell,k}\bigg{[}\prod_{j=1}^{s}(Q_{j})_{\ell,k}\bigg{]}^{-1}\Gamma^{-1}_{k}(\ell B+C)\\
&=\sum_{\ell=0}^{\infty}(A+k\ell I)\mathbf{\Psi}_{\ell,k}(z).\label{2.6}
\end{split}
\end{equation}
%-----------------------------------------------------------------------------------------------------------------------------------------------------------------
where $\mathbf{\Psi}_{\ell,k}(z)=\frac{z^{\ell}}{\ell !}(A)_{\ell,k}\prod_{i=1}^{r}(P_{i})_{\ell,k}\bigg{[}\prod_{j=1}^{s}(Q_{j})_{\ell,k}\bigg{]}^{-1}\Gamma^{-1}_{k}(\ell B+C)$

%-----------------------------------------------------------------------------------------------------------------------------------------------------------------
Similarly, we get
%-----------------------------------------------------------------------------------------------------------------------------------------------------------------
\begin{equation}
\begin{split}
&\;_{r+1}R_{s,k}(A-kI)=(A-kI)\sum_{\ell=0}^{\infty}\bigg{(}A+k(\ell-1)I\bigg{)}^{-1}\mathbf{\Psi}_{\ell,k}(z),\\
&P_{i}\;_{r+1}R_{s,k}(P_{i}+kI)=\sum_{\ell=0}^{\infty}(P_{i}+k\ell I)\mathbf{\Psi}_{\ell,k}(z),\\
&\;_{r+1}R_{s,k}(P_{i}-kI)=(P_{i}-kI)\sum_{\ell=0}^{\infty}\bigg{(}P_{i}+k(\ell-1)I\bigg{)}^{-1}W_{\ell,k}(z)\mathbf{\Psi}_{\ell,k}(z),\\
&\;_{r+1}R_{s,k}(Q_{j}+kI)=Q_{j}\sum_{\ell=0}^{\infty}\bigg{(}Q_{j}+k\ell I\bigg{)}^{-1}\mathbf{\Psi}_{\ell,k}(z),\\
&\big{(}Q_{j}-kI\big{)}\;_{r+1}R_{s,k}(Q_{j}-kI)=\sum_{\ell=0}^{\infty}(Q_{j}+k(\ell-1)I)\mathbf{\Psi}_{\ell,k}(z).\label{2.7}
\end{split}
\end{equation}
%-----------------------------------------------------------------------------------------------------------------------------------------------------------------
For all integers $n\ge 1$, we have
%-----------------------------------------------------------------------------------------------------------------------------------------------------------------
\begin{equation}
\begin{split}
\;_{r+1}R_{s,k}(A+nkI)&=\sum_{\ell=0}^{\infty}\prod_{\mu=1}^{n}\bigg{(}A+k(\mu-1)I\bigg{)}^{-1}(A+k(\ell+\mu-1)I)\mathbf{\Psi}_{\ell,k}(z),\\
\;_{r+1}R_{s,k}(A-nkI)&=\sum_{\ell=0}^{\infty}\prod_{\mu=1}^{n}(A-k\mu I)\bigg{(}A+k(\ell-\mu)I\bigg{)}^{-1}\mathbf{\Psi}_{\ell,k}(z),\\
\;_{r+1}R_{s,k}(P_{i}+nkI)&=\sum_{\ell=0}^{\infty}\prod_{\mu=1}^{n}\bigg{(}P_{i}+k(\mu-1)I\bigg{)}^{-1}(P_{i}+k(\ell+\mu-1)I)\mathbf{\Psi}_{\ell,k}(z),\\
\;_{r+1}R_{s,k}(P_{i}-nkI)&=\sum_{\ell=0}^{\infty}\prod_{\mu=1}^{n}(P_{i}-k\mu I)\bigg{(}P_{i}+k(\ell-\mu)I\bigg{)}^{-1}\mathbf{\Psi}_{\ell,k}(z),\\
\;_{r+1}R_{s,k}(Q_{j}+nkI)&=\sum_{\ell=0}^{\infty}\prod_{\mu=1}^{n}(Q_{j}+k(\mu-1)I)\bigg{(}Q_{j}+k(\ell+\mu-1)I\bigg{)}^{-1}\mathbf{\Psi}_{\ell,k}(z),\\
\;_{r+1}R_{s,k}(Q_{j}-nkI)&=\sum_{\ell=0}^{\infty}\prod_{\mu=1}^{n}\bigg{(}Q_{j}-k\mu I\bigg{)}^{-1}(Q_{j}+k(\ell-\mu)I)\mathbf{\Psi}_{\ell,k}(z).\label{2.8}
\end{split}
\end{equation}
%-----------------------------------------------------------------------------------------------------------------------------------------------------------------
Using the differential operator $\theta=z\frac{d}{dz}$, we get
%-----------------------------------------------------------------------------------------------------------------------------------------------------------------
\begin{equation}
\begin{split}
(k\theta I+A)\;_{r+1}R_{s,k}&=\sum_{\ell=0}^{\infty}(A+k\ell I)\mathbf{\Psi}_{\ell,k}(z),\\
(k\theta I+P_{i})\;_{r+1}R_{s,k}&=\sum_{\ell=0}^{\infty}(P_{i}+k\ell I)\mathbf{\Psi}_{\ell,k}(z).\label{2.9}
\end{split}
\end{equation}
%-----------------------------------------------------------------------------------------------------------------------------------------------------------------
From (\ref{2.6}) and (\ref{2.9}), we have
%-----------------------------------------------------------------------------------------------------------------------------------------------------------------
\begin{equation}
\begin{split}
&(k\theta I+A)\;_{r+1}R_{s,k}=A\;_{r+1}R_{s,k}(A+k I),\\
&(k\theta I+P_{i})\;_{r+1}R_{s,k}=P_{i}\;_{r+1}R_{s,k}(P_{i}+k I)\,;\,i=1,2,\ldots,r.\label{2.10}
\end{split}
\end{equation}
%-----------------------------------------------------------------------------------------------------------------------------------------------------------------
Similarly, we get
%-----------------------------------------------------------------------------------------------------------------------------------------------------------------
\begin{equation}
\begin{split}
(k\theta I+Q_{j}-k I)\;_{r+1}R_{s,k}=(Q_{j}-kI)\;_{r+1}R_{s,k}(Q_{j}-kI)\,;\,j=1,2,\ldots,s.\label{2.11}
\end{split}
\end{equation}
%-----------------------------------------------------------------------------------------------------8------------------------------------------------------------
The elimination of $\theta\;_{r+1}R_{s,k}$ from (\ref{2.10}) and (\ref{2.11}), we obtain the recurrence matrix relations (\ref{2.3}), (\ref{2.4}) and (\ref{2.5})
%---------------------------------------------------------------------------------------------------------------------------------------------------------------------------------------------
\end{proof}
%---------------------------------------------------------------------------------------------------------------------------------------------------------------------------------------------
%---------------------------------------------------------------------------------------------------------------------------------------------------------------------------------------------
\begin{theorem} The following differential formulas hold true for $\;_{r+1}R_{s,k}$
%---------------------------------------------------------------------------------------------------------------------------------------------------------------------------------------------
\begin{eqnarray}
\begin{split}
\;_{r+1}R_{s,k}=C\;_{r+1}R_{s,k}(C+kI)+zB\frac{d}{dz}\;_{r+1}R_{s,k}(C+kI)\label{2.12}
\end{split}
\end{eqnarray}
%---------------------------------------------------------------------------------------------------------------------------------------------------------------------------------------------
and
%-----------------------------------------------------------------------------------------------------------------------------------------------------------------
\begin{equation}
\begin{split}
D_{z}^{\mu}\;_{r+1}R_{s,k}=&(A)_{\mu,k}\Pi_{i=1}^{r}(P_{i})_{\mu,k}\bigg{[}\Pi_{j=1}^{s}(Q_{j})_{\mu,k}\bigg{]}^{-1}\;_{r+1}R_{s,k}(A+\mu kI,P_{1}+\mu kI,\ldots,P_{r}+\mu kI;\\
&Q_{1}+\mu kI,\ldots,Q_{s}+\mu kI;B,\mu B+C,z).\label{2.13}
\end{split}
\end{equation}
%---------------------------------------------------------------------------------------------------------------------------------------------------------------------------------------------

%---------------------------------------------------------------------------------------------------------------------------------------------------------------------------------------------
\end{theorem}
%---------------------------------------------------------------------------------------------------------------------------------------------------------------------------------------------
\begin{proof}
%---------------------------------------------------------------------------------------------------------------------------------------------------------------------------------------------
The right hand side in (\ref{2.12}) and using (\ref{2.1}), we get
%-----------------------------------------------------------------------------------------------------------------------------------------------------------------
\begin{equation*}
\begin{split}
&C\;_{r+1}R_{s,k}(C+kI)+zB\frac{d}{dz}\;_{r+1}R_{s,k}(C+kI)\\
&=C\;_{r+1}R_{s,k}(C+kI)+zB\bigg{[}\sum_{\ell=0}^{\infty}\frac{\ell z^{\ell-1}}{\ell !}(A)_{\ell,k}\prod_{i=1}^{r}(P_{i})_{\ell,k}\bigg{[}\prod_{j=1}^{s}(Q_{j})_{\ell,k}\bigg{]}^{-1}\Gamma_{k}^{-1}(\ell B+C+kI)\bigg{]}\\
&=C\;_{r+1}R_{s,k}(C+kI)+\sum_{\ell=0}^{\infty}\frac{z^{\ell}}{\ell !}(A)_{\ell,k}\prod_{i=1}^{r}(P_{i})_{\ell,k}\bigg{[}\prod_{j=1}^{s,k}(Q_{j})_{\ell,k}\bigg{]}^{-1}\Gamma_{k}^{-1}(\ell B+C)\\
&-C\sum_{\ell=0}^{\infty}\frac{ z^{\ell}}{\ell !}(A)_{\ell,k}\prod_{i=1}^{r}(P_{i})_{\ell,k}\bigg{[}\prod_{j=1}^{s}(Q_{j})_{\ell,k}\bigg{]}^{-1}\Gamma_{k}^{-1}(\ell B+C+kI)\bigg{]}=\;_{r+1}R_{s,k}.
\end{split}
\end{equation*}
%---------------------------------------------------------------------------------------------------------------------------------------------------------------------------------------
%---------------------------------------------------------------------------------------------------------------------------------------------------------------------------------------------
Differentiating (\ref{2.1}) with respect to $z$, we get
%-----------------------------------------------------------------------------------------------------------------------------------------------------------------
\begin{equation*}
\begin{split}
D_{z}\;_{r+1}R_{s,k}=&A\Pi_{i=1}^{r}P_{i}\bigg{(}\Pi_{j=1}^{s}Q_{j}\bigg{)}^{-1}\;_{r+1}R_{s,k}(A+kI,P_{1}+kI,\\
&\ldots,P_{r}+kI;Q_{1}+kI,\ldots,Q_{s}+kI;B,B+C,z)\,;\, D_{z}=\frac{d}{dz}.
\end{split}
\end{equation*}
%-----------------------------------------------------------------------------------------------------------------------------------------------------------------
By repeating the above process that $\mu$ times, we get (\ref{2.13}).
%-----------------------------------------------------------------------------------------------------------------------------------------------------------------
%---------------------------------------------------------------------------------------------------------------------------------------------------------------------------------------------
\end{proof}
%---------------------------------------------------------------------------------------------------------------------------------------------------------------------------------------------
%-----------------------------------------------------------------------------------------------------------------------------------------------------------------------------------
%-----------------------------------------------------------------------------------------------------------------------------------------------------------------------------------
\begin{remark}
%-----------------------------------------------------------------------------------------------------------------------------------------------------------------------------------
\begin{enumerate}
  \item If $k=1$, $A=B=C=I$ in (\ref{2.3})-(\ref{2.5}), we get the results for the generalized hypergeometric matrix functions \cite{sa1}.
  \item If $k=1$ in (\ref{2.12})-(\ref{2.13}), we get the contiguous function relations for the $\;_{r+1}R_{s}$ matrix function \cite{sd2, skc}.
\end{enumerate}
%-----------------------------------------------------------------------------------------------------------------------------------------------------------------------------------
\end{remark}
%---------------------------------------------------------------------------------------------------------------------------------------------------------------------------------------------
%---------------------------------------------------------------------------------------------------------------------------------------------------------------------------------------------
\begin{theorem}
%---------------------------------------------------------------------------------------------------------------------------------------------------------------------------------------------
The $\;_{r+1}R_{s,k}$ matrix function has the following differential properties:
%---------------------------------------------------------------------------------------------------------------------------------------------------------------------------------------------
\begin{equation}
\begin{split}
&\bigg{(}\frac{d}{dz}\bigg{)}^{\mu}\bigg{[}z^{\frac{C}{k}-I}\;_{r+1}R_{s,k}(A,P_{1},P_{2},\ldots,P_{r};Q_{1},Q_{2},\ldots,Q_{s};B,C;z^{\frac{B}{k}})\bigg{]}\\
&=\frac{1}{k^{\mu}}z^{\frac{C}{k}-(\mu+1)I}\;_{r+1}R_{s,k}(A,P_{1},P_{2},\ldots,P_{r};Q_{1},Q_{2},\ldots,Q_{s};B,C-\mu I;z^{\frac{B}{k}}),\label{2.14}
\end{split}
\end{equation}
%---------------------------------------------------------------------------------------------------------------------------------------------------------------------------------------------
%---------------------------------------------------------------------------------------------------------------------------------------------------------------------------------------------
\begin{equation}
\begin{split}
&\bigg{(}\frac{d}{dz}\bigg{)}^{\mu}\bigg{[}z^{\frac{A}{k}+\mu I-I}\;_{r+1}R_{s,k}(A,P_{1},P_{2},\ldots,P_{r};Q_{1},Q_{2},\ldots,Q_{s};B,C;z)\bigg{]}\\
&=\frac{1}{k^{\mu}}(A)_{\mu,k}z^{\frac{A}{k}-I}\;_{r+1}R_{s,k}(A+\mu kI,P_{1},P_{2},\ldots,P_{r};Q_{1},Q_{2},\ldots,Q_{s};B,C;z),\label{2.15}
\end{split}
\end{equation}
%---------------------------------------------------------------------------------------------------------------------------------------------------------------------------------------------
%---------------------------------------------------------------------------------------------------------------------------------------------------------------------------------------------
\begin{equation}
\begin{split}
&\bigg{(}\frac{d}{dz}\bigg{)}^{\mu}\bigg{[}z^{\frac{P_{i}}{k}+\mu I-I}\;_{r+1}R_{s,k}(A,P_{1},P_{2},\ldots,P_{r};Q_{1},Q_{2},\ldots,Q_{s};B,C;z)\bigg{]}\\
&=\frac{1}{k^{\mu}}(P_{i})_{\mu,k}z^{\frac{P_{i}}{k}-I}\;_{r+1}R_{s,k}(P_{i}+\mu kI),i=1,2,\ldots r\label{2.16}
\end{split}
\end{equation}
%---------------------------------------------------------------------------------------------------------------------------------------------------------------------------------------------
and
%---------------------------------------------------------------------------------------------------------------------------------------------------------------------------------------------
\begin{equation}
\begin{split}
&\bigg{(}\frac{d}{dz}\bigg{)}^{\mu}\bigg{[}z^{\frac{Q_{j}}{k}}\;_{r+1}R_{s,k}(Q_{j}+kI;z)\bigg{]}\\
&=\frac{1}{k^{\mu}}\Gamma_{k}(Q_{j})\Gamma_{k}^{-1}(Q_{j}-\mu k I)z^{\frac{Q_{j}}{k}-\mu I}\;_{r+1}R_{s,k}(Q_{j}-(\mu-1)k I;z),1\leq j\leq s.\label{2.17}
\end{split}
\end{equation}
%---------------------------------------------------------------------------------------------------------------------------------------------------------------------------------------------
%---------------------------------------------------------------------------------------------------------------------------------------------------------------------------------------------
\end{theorem}
%---------------------------------------------------------------------------------------------------------------------------------------------------------------------------------------------
%---------------------------------------------------------------------------------------------------------------------------------------------------------------------------------------------
\begin{proof}
%-----------------------------------------------------------------------------------------------------------------------------------------------------------------
%-----------------------------------------------------------------------------------------------------------------------------------------------------------------
To prove (\ref{2.14}). Multiplying the equation (\ref{2.1}) by $z^{\frac{C}{k}-I}$ and differentiating with respect to $z$, we get
%---------------------------------------------------------------------------------------------------------------------------------------------------------------------------------------------
\begin{equation*}
\begin{split}
&\bigg{(}\frac{d}{dz}\bigg{)}\bigg{[}z^{\frac{C}{k}-I}\;_{r+1}R_{s,k}(A,P_{1},P_{2},\ldots,P_{r};Q_{1},Q_{2},\ldots,Q_{s};B,C;z^{\frac{B}{k}})\bigg{]}\\
&=\frac{1}{k}z^{\frac{C}{k}-2I}\;_{r+1}R_{s,k}(A,P_{1},P_{2},\ldots,P_{r};Q_{1},Q_{2},\ldots,Q_{s};B,C-I;z^{\frac{B}{k}}).
\end{split}
\end{equation*}
%-----------------------------------------------------------------------------------------------------------------------------------------------------------------
By repeating the above relation $\mu$ times, we get (\ref{2.14}). On the same parallel lines, we get the results (\ref{2.15})-(\ref{2.17}). So the detailed account of proof is omitted.
%---------------------------------------------------------------------------------------------------------------------------------------------------------------------------------------------
\end{proof}
%---------------------------------------------------------------------------------------------------------------------------------------------------------------------------------------------

%---------------------------------------------------------------------------------------------------------------------------------------------------------------------------------------------
\begin{theorem}
%---------------------------------------------------------------------------------------------------------------------------------------------------------------------------------------------
The $\;_{r+1}R_{s,k}$ matrix functions satisfies the differential equation
%---------------------------------------------------------------------------------------------------------------------------------------------------------------------------------------------
\begin{equation}
\begin{split}
\theta\; \prod_{j=1}^{s}(k\theta \; I+Q_{j}-k I)\;_{r+1}R_{s,k}-z(k\theta I+A)\prod_{i=1}^{r}(k\theta\; I+P_{i})\;_{r+1}R_{s,k}(B+C)=\textbf{0},\label{2.18}
\end{split}
\end{equation}
%---------------------------------------------------------------------------------------------------------------------------------------------------------------------------------------------
where $\textbf{0}$ is the null matrix in $\Bbb{C}^{N\times N}$.
%---------------------------------------------------------------------------------------------------------------------------------------------------------------------------------------------
\end{theorem}
%---------------------------------------------------------------------------------------------------------------------------------------------------------------------------------------------
%---------------------------------------------------------------------------------------------------------------------------------------------------------------------------------------------
\begin{proof}
%---------------------------------------------------------------------------------------------------------------------------------------------------------------------------------------------
%-----------------------------------------------------------------------------------------------------------------------------------------------------------------
Using the Euler differential operator $\theta=z\frac{d}{dz}$, we get
%-----------------------------------------------------------------------------------------------------------------------------------------------------------------
\begin{equation*}
\begin{split}
\theta\; \prod_{j=1}^{s}(k\theta\; I+Q_{j}-kI)\;_{r+1}R_{s,k}&=\sum_{\ell=1}^{\infty}\frac{\ell\;z^{\ell}}{\ell !}\prod_{j=1}^{s}(k\ell I+Q_{j}-k I)(A)_{\ell,k}\prod_{i=1}^{r}(P_{i})_{\ell,k}\bigg{(}\prod_{j=1}^{s}(Q_{j})_{\ell,k}\bigg{]}^{-1}\Gamma^{-1}_{k}(\ell B+C)\\
&=\sum_{\ell=1}^{\infty}\frac{\;z^{\ell}}{(\ell-1)!}(A)_{\ell,k}\prod_{i=1}^{r}(P_{i})_{\ell,k}\bigg{(}\prod_{j=1}^{s}(Q_{j})_{\ell-1,k}\bigg{)}^{-1}\Gamma^{-1}_{k}(\ell B+C).
\end{split}
\end{equation*}
%-----------------------------------------------------------------------------------------------------------------------------------------------------------------
Replace $\ell$ by $\ell+1$, we have
%-----------------------------------------------------------------------------------------------------------------------------------------------------------------
\begin{equation*}
\begin{split}
\theta\; \prod_{j=1}^{s}(k\theta\; I+Q_{j}-k I)\;_{r+1}R_{s,k}&=\sum_{\ell=0}^{\infty}\frac{\;z^{\ell+1}}{\ell !}(A)_{\ell+1,k}\prod_{i=1}^{r}(P_{i})_{\ell+1,k}\bigg{(}\prod_{j=1}^{s}(Q_{j})_{\ell,k}\bigg{)}^{-1}\Gamma^{-1}_{k}(\ell B+B+C)\\
&=z(k\theta I+A)\prod_{i=1}^{r}(k\theta\; I+P_{i})\;_{r+1}R_{s,k}(B+C).
\end{split}
\end{equation*}
%---------------------------------------------------------------------------------------------------------------------------------------------------------------------------------------------
\end{proof}
%-----------------------------------------------------------------------------------------------------------------------------------------------------------------------------------
%-----------------------------------------------------------------------------------------------------------------------------------------------------------------------------------
\begin{theorem}
%-----------------------------------------------------------------------------------------------------------------------------------------------------------------------------------
The following integral representations for $\;_{r+1}R_{s,k}$ matrix function hold true:
%-----------------------------------------------------------------------------------------------------------------------------------------------------------------------------------
\begin{equation}
\begin{split}
&\;_{r+1}R_{s,k}\left(A,P_{1},P_{2},\ldots ,P_{r};Q_{1},Q_{2},\ldots ,Q_{s};B,C,z\right)\\
&=\frac{1}{k}\Gamma_{k} \left( Q_{j}\right)\Gamma_{k} ^{-1}\left(A\right) \Gamma_{k} ^{-1}\left( Q_{j}-A\right)
\int_{0}^{1}\xi^{\frac{1}{k}A-I}\left( 1-\xi\right)
^{\frac{1}{2}(Q_{j}-A)-I} \\
&\times\;_{r}R_{s-1,k}\left(
\begin{array}{c}
P_{1},\ldots ,P_{i-1},P_{i},P_{i+1}\ldots ,P_{r}; \\
Q_{1},\ldots ,Q_{j-1},Q_{j+1}\ldots ,Q_{s}
\end{array}
;B,C,z\xi\right)d\xi,\label{2.19}
\end{split}
\end{equation}
%-----------------------------------------------------------------------------------------------------------------------------------------------------------------------------------
\begin{equation}
\begin{split}
&\;_{r+1}R_{s,k}\left(A,P_{1},P_{2},\ldots ,P_{r};Q_{1},Q_{2},\ldots ,Q_{s};B,C,z\right)\\
&=\frac{1}{k}\Gamma_{k} \left( Q_{j}\right)\Gamma_{k} ^{-1}\left( P_{i}\right) \Gamma_{k} ^{-1}\left( Q_{j}-P_{i}\right)
\int_{0}^{1}\xi^{\frac{1}{k}P_{i}-I}\left( 1-\xi\right)
^{\frac{1}{2}(Q_{j}-P_{i})-I} \\
&\times\;_{r}R_{s-1,k}\left(
\begin{array}{c}
A,P_{1},\ldots ,P_{i-1},P_{i+1}\ldots ,P_{r}; \\
Q_{1},\ldots ,Q_{j-1},Q_{j+1}\ldots ,Q_{s}
\end{array}
;B,C,z\xi\right)d\xi,\label{2.20}
\end{split}
\end{equation}
%-----------------------------------------------------------------------------------------------------------------------------------------------------------------------------------
%---------------------------------------------------------------------------------------------------------------------------------------------------------------------------------------------
\begin{equation}
\begin{split}
&\;_{r+1}R_{s,k}(A,P_{1},P_{2},\ldots,P_{r};Q_{1},Q_{2},\ldots,Q_{s};P,Q,z)=\Gamma_{k}^{-1}(A)\\
&\int_{0}^{\infty}\xi^{A-I}e^{-\frac{\xi^{k}}{k}}\;_{r}R_{s,k}(P_{1},P_{2},\ldots,P_{r};Q_{1},Q_{2},\ldots,Q_{s};B,C,z\xi^{k})d\xi\label{2.21}
\end{split}
\end{equation}
%---------------------------------------------------------------------------------------------------------------------------------------------------------------------------------------------
and
%---------------------------------------------------------------------------------------------------------------------------------------------------------------------------------------------
\begin{equation}
\begin{split}
&\;_{r+1}R_{s,k}(A,P_{1},P_{2},\ldots,P_{r};Q_{1},Q_{2},\ldots,Q_{s};P,Q,z)=\Gamma_{k}^{-1}(P_{i})\\
&\int_{0}^{\infty}\xi^{P_{i}-I}e^{-\frac{\xi^{k}}{k}}\;_{r}R_{s,k}(A,P_{1},\ldots ,P_{i-1},P_{i+1}\ldots ,P_{r};Q_{1},Q_{2},\ldots,Q_{s};B,C,z\xi^{k})d\xi.\label{2.22}
\end{split}
\end{equation}
%---------------------------------------------------------------------------------------------------------------------------------------------------------------------------------------------
\end{theorem}
%-----------------------------------------------------------------------------------------------------------------------------------------------------------------------------------
\begin{proof}
%-----------------------------------------------------------------------------------------------------------------------------------------------------------------------------------
By using (\ref{1.3}) and (\ref{1.5}), we get
%---------------------------------------------------------------------------------------------------------------------------------------------------------------------------------------------
\begin{equation}
\begin{split}
&(A)_{\ell,k}\bigg{[}(Q_{j})_{\ell,k}\bigg{]}^{-1}=\Gamma_{k}(Q_{j})\Gamma_{k}(A+k\ell)\Gamma_{k}^{-1}(A)\Gamma_{k}^{-1}(Q_{j}+k\ell)\\
&=\Gamma_{k}(Q_{j})\Gamma_{k}^{-1}(A)\Gamma_{k}^{-1}(Q_{j}-A)\Bbb{B}_{k}(A+k\ell I,Q_{j}-A)\\
&=\frac{1}{k}\Gamma_{k}(Q_{j})\Gamma_{k}^{-1}(A)\Gamma_{k}^{-1}(Q_{j}-A)\int_{0}^{1}\xi^{\frac{1}{k}A+\left( \ell-1\right) I}\left( 1-\xi\right)
^{\frac{1}{k}(Q_{j}-A)-I}d\xi.\label{2.23}
\end{split}
\end{equation}
%-----------------------------------------------------------------------------------------------------------------------------------------------------------------------------------
Using the above equation (\ref{2.23}) and (\ref{2.1}), we obtain (\ref{2.19}). In a similar way, we obtain to the desired results (\ref{2.20})-(\ref{2.22}).
%-----------------------------------------------------------------------------------------------------------------------------------------------------------------------------------
%-----------------------------------------------------------------------------------------------------------------------------------------------------------------------------------
\end{proof}

%---------------------------------------------------------------------------------------------------------------------------------------------------------------------------------------------
\begin{theorem} The following integrals representations hold true:
%---------------------------------------------------------------------------------------------------------------------------------------------------------------------------------------------
%---------------------------------------------------------------------------------------------------------------------------------------------------------------------------------------------
\begin{equation}
\begin{split}
&\int_{0}^{x}\xi^{\frac{Q_{j}}{k}-I}(x-\xi)^{\frac{E}{k}-I}\;_{r+1}R_{s,k}(A,P_{1},P_{2},\ldots,P_{r};Q_{1},Q_{2},\ldots,Q_{s};B,C,z\xi)d\xi\\
&=k\Bbb{B}_{k}(Q_{j},E)x^{\frac{Q_{j}+E}{k}-I}\;_{r+1}R_{s,k}(A,P_{1},P_{2},\ldots,P_{r};Q_{j}+E;B,C;zx),1\leq j\leq s\label{2.24}
\end{split}
\end{equation}
%---------------------------------------------------------------------------------------------------------------------------------------------------------------------------------------------
and
%---------------------------------------------------------------------------------------------------------------------------------------------------------------------------------------------
\begin{equation}
\begin{split}
&\int_{a}^{x}(x-\xi)^{\frac{E}{k}-I}(\xi-a)^{\frac{Q_{j}}{k}-I}\;_{r+1}R_{s,k}(A,P_{1},P_{2},\ldots,P_{r};Q_{j};B,C,a-\xi)d\xi\\
&=k\Bbb{B}_{k}(E,Q_{j})(x-a)^{\frac{Q_{j}+E}{k}-I}\;_{r+1}R_{s,k}(A,P_{1},P_{2},\ldots,P_{r};Q_{j}+E;B,C;a-x),1\leq j\leq s.\label{2.25}
\end{split}
\end{equation}
%---------------------------------------------------------------------------------------------------------------------------------------------------------------------------------------------
%---------------------------------------------------------------------------------------------------------------------------------------------------------------------------------------------
\end{theorem}
%---------------------------------------------------------------------------------------------------------------------------------------------------------------------------------------------
%---------------------------------------------------------------------------------------------------------------------------------------------------------------------------------------------
\begin{proof}
%---------------------------------------------------------------------------------------------------------------------------------------------------------------------------------------------
Taking left hand side of (\ref{2.24}), we get
%---------------------------------------------------------------------------------------------------------------------------------------------------------------------------------------------
\begin{equation*}
\begin{split}
&\int_{0}^{x}\xi^{\frac{Q_{j}}{k}-I}(x-\xi)^{\frac{E}{k}-I}\sum_{\ell=0}^{\infty}\frac{(z\xi)^{\ell}}{\ell !}(A)_{\ell,k}\prod_{i=1}^{r}(P_{i})_{\ell,k}\bigg{[}\prod_{j=2}^{s}(Q_{j})_{\ell,k}\bigg{]}^{-1}\Gamma^{-1}_{k}(\ell B+C)d\xi\\
&=\sum_{\ell=0}^{\infty}\int_{0}^{x}\xi^{\frac{Q_{j}}{k}+(\ell-1)}(x-\xi)^{\frac{E}{k}-I}\frac{z^{\ell}}{\ell !}(A)_{\ell,k}\prod_{i=1}^{r}(P_{i})_{\ell,k}\bigg{[}\prod_{j=2}^{s}(Q_{j})_{\ell,k}\bigg{]}^{-1}\Gamma^{-1}_{k}(\ell B+C)d\xi.
\end{split}
\end{equation*}
%---------------------------------------------------------------------------------------------------------------------------------------------------------------------------------------------
Changing order of summation and integration with put $\xi=xu$ and $d\xi=xdu$, we get
%---------------------------------------------------------------------------------------------------------------------------------------------------------------------------------------------
\begin{equation*}
\begin{split}
&\sum_{\ell=0}^{\infty}x^{\frac{Q_{j}+E}{k}+(\ell-1)I}\int_{0}^{1}u^{\frac{Q_{j}}{k}+\ell-I}(1-u)^{\frac{E}{k}-I}\frac{z^{\ell}}{\ell !}(A)_{\ell,k}\prod_{i=1}^{r}(P_{i})_{\ell,k}\bigg{[}\prod_{j=2}^{s}(Q_{j})_{\ell,k}\bigg{]}^{-1}\Gamma^{-1}_{k}(\ell B+C)du\\
&=\sum_{\ell=0}^{\infty}x^{\frac{Q_{j}+E}{k}+(\ell-1)I}\frac{z^{\ell}}{\ell !}(A)_{\ell,k}\prod_{i=1}^{r}(P_{i})_{\ell,k}\bigg{[}\prod_{j=2}^{s}(Q_{j})_{\ell,k}\bigg{]}^{-1}\Gamma^{-1}_{k}(\ell B+C)k\Bbb{B}_{k}(Q_{j}+\ell kI,E)\\
&=k\sum_{\ell=0}^{\infty}x^{\frac{Q_{j}+E}{k}+(\ell-1)I}\frac{z^{\ell}}{\ell !}(A)_{\ell,k}\prod_{i=1}^{r}(P_{i})_{\ell,k}\bigg{[}\prod_{j=2}^{s}(Q_{j})_{\ell,k}\bigg{]}^{-1}\Gamma^{-1}_{k}(\ell B+C)\Gamma_{k}(Q_{j}+\ell kI)\Gamma_{k}(E)\Gamma_{k}^{-1}(Q_{j}+E+\ell kI)
\end{split}
\end{equation*}
%---------------------------------------------------------------------------------------------------------------------------------------------------------------------------------------------
%---------------------------------------------------------------------------------------------------------------------------------------------------------------------------------------------
Which is desired result (\ref{2.24}). Letting $t=\frac{\xi-a}{x-a}$,$\xi-a=(x-a)t$,$x-\xi=(x-a)(1-t)$ and using the $k$-beta matrix function (\ref{1.3}) in the left hand side of (\ref{2.25}), we obtain (\ref{2.25}).
%---------------------------------------------------------------------------------------------------------------------------------------------------------------------------------------------
\end{proof}
%---------------------------------------------------------------------------------------------------------------------------------------------------------------------------------------------
%---------------------------------------------------------------------------------------------------------------------------------------------------------------------------------------------
\begin{theorem} Let $\mu$ be a positive integer, then the $\;_{r+1}R_{s,k}$ matrix function satisfies the following Euler-Type integral representation
%-----------------------------------------------------------------------------------------------------------------------------------------------------------------------------------
\begin{equation}
\begin{split}
&\;_{r+\mu}R_{s+\mu}\left(A,P_{1},P_{2},\ldots ,P_{r},\Delta(E;\mu);Q_{1},Q_{2},\ldots ,Q_{s},\triangle(E+M;\mu);B,C,cz^{\frac{\mu}{k}}\right)\\
&=\frac{1}{k}z^{I-\frac{E+M}{k}}\Gamma ^{-1}_{k}\left( E\right) \Gamma_{k}\left(E+M\right)
\Gamma^{-1}_{k}\left( M\right) \int_{0}^{z}\xi^{\frac{E}{k}-I}\left( z-\xi\right)
^{\frac{M}{k}-I} \\
&\times\;_{r+1}R_{s,k}\left(
\begin{array}{c}
A,P_{1},P_{2},\ldots ,P_{r}; \\
Q_{1},Q_{2},\ldots ,Q_{s}
\end{array}
;B,C,c\xi^{\frac{\mu}{k}}\right)d\xi,\label{2.26}
\end{split}
\end{equation}
%---------------------------------------------------------------------------------------------------------------------------------------------------------------------------------------------
where $\Delta(E,\mu)$ is given the array of a parameters
%---------------------------------------------------------------------------------------------------------------------------------------------------------------------------------------------
\begin{equation*}
\begin{split}
\Delta(E,\mu)=\frac{1}{\mu}E,\frac{1}{\mu}(E+kI),\frac{1}{\mu}(E+2kI),\ldots,\frac{1}{\mu}(E+k(\mu-1)I).
\end{split}
\end{equation*}
%---------------------------------------------------------------------------------------------------------------------------------------------------------------------------------------------
and
%-----------------------------------------------------------------------------------------------------------------------------------------------------------------------------------
\begin{equation}
\begin{split}
&\;_{r+\mu+\imath}R_{s+\mu+\imath}\left(A,P_{1},P_{2},\ldots ,P_{r},\Delta(E;\mu),\Delta(M;\imath);Q_{1},Q_{2},\ldots ,Q_{s},\triangle(E+M;\mu+\imath);B,C,\frac{\alpha\mu^{\mu}\imath^{\imath}}{(\mu+\imath)^{\mu+\imath}}\right)\\
&=\Gamma ^{-1}_{k}\left( E\right) \Gamma_{k}\left(E+M\right)
\Gamma^{-1}_{k}\left( M\right) \int_{0}^{1}\xi^{\frac{E}{k}-I}\left(1-\xi\right)
^{\frac{M}{k}-I} \\
&\times\;_{r+1}R_{s,k}\left(
\begin{array}{c}
A,P_{1},P_{2},\ldots ,P_{r}; \\
Q_{1},Q_{2},\ldots ,Q_{s}
\end{array}
;B,C,\alpha \xi^{\frac{\mu}{k}}(1-t)^{\frac{\imath}{k}}\right)d\xi.\label{2.27}
\end{split}
\end{equation}
%---------------------------------------------------------------------------------------------------------------------------------------------------------------------------------------------
\end{theorem}
%---------------------------------------------------------------------------------------------------------------------------------------------------------------------------------------------
\begin{proof}
%---------------------------------------------------------------------------------------------------------------------------------------------------------------------------------------------
If we apply the substitution $\chi=zu$, $\chi=z\xi$, $\chi=zd\xi$ and use the $k$-Beta matrix function, we obtain our
desired results (\ref{2.26})-(\ref{2.27}).
%---------------------------------------------------------------------------------------------------------------------------------------------------------------------------------------------
\end{proof}
%---------------------------------------------------------------------------------------------------------------------------------------------------------------------------------------------
\begin{theorem}
%---------------------------------------------------------------------------------------------------------------------------------------------------------------------------------------------
The following results hold true:
%---------------------------------------------------------------------------------------------------------------------------------------------------------------------------------------------
%---------------------------------------------------------------------------------------------------------------------------------------------------------------------------------------------
\begin{equation}
\begin{split}
&\mathfrak{B}\bigg{[}\;_{r+1}R_{s,k}(A+E,P_{1},P_{2},\ldots,P_{r};Q_{1},Q_{2},\ldots,Q_{s};B,C,zt);A,E\bigg{]}\\
&=\frac{1}{k}\int_{0}^{1}t^{\frac{A}{k}-I}(1-t)^{\frac{E}{k}-I}\;_{r+1}R_{s,k}(A+E,P_{1},P_{2},\ldots,P_{r};Q_{1},Q_{2},\ldots,Q_{s};B,C,zt)dt\\
&=\Gamma_{k}(A)\Gamma_{k}(E)\Gamma_{k}^{-1}(A+E)\;_{r+1}R_{s,k}(A,P_{1},P_{2},\ldots,P_{r};Q_{1},Q_{2},\ldots,Q_{s};B,C;z)\label{2.28}
\end{split}
\end{equation}
%---------------------------------------------------------------------------------------------------------------------------------------------------------------------------------------------
and
%---------------------------------------------------------------------------------------------------------------------------------------------------------------------------------------------
\begin{equation}
\begin{split}
&\mathfrak{B}\bigg{[}\;_{r+1}R_{s,k}(A,P_{1},P_{2},\ldots, P_{i-1},P_{i}+E,P_{i+1},\ldots,P_{r};Q_{1},Q_{2},\ldots,Q_{s};B,C,zt);P_{i},E\bigg{]}\\
&=\frac{1}{k}\int_{0}^{1}t^{\frac{P_{i}}{k}-I}(1-t)^{\frac{E}{k}-I}\;_{r+1}R_{s,k}(A,P_{1},P_{2},\ldots, P_{i-1},P_{i}+E,P_{i+1},,\ldots,P_{r};Q_{1},Q_{2},\ldots,Q_{s};B,C,zt)dt\\
&=\Gamma_{k}(P_{i})\Gamma_{k}(E)\Gamma_{k}^{-1}(P_{i}+E)\;_{r+1}R_{s,k}(A,P_{1},P_{2},\ldots, P_{i-1},P_{i+1},\ldots,P_{r};Q_{1},Q_{2},\ldots,Q_{s};B,C;z),\label{2.29}
\end{split}
\end{equation}
%---------------------------------------------------------------------------------------------------------------------------------------------------------------------------------------------
\end{theorem}
%---------------------------------------------------------------------------------------------------------------------------------------------------------------------------------------------
%---------------------------------------------------------------------------------------------------------------------------------------------------------------------------------------------
\begin{proof}
%---------------------------------------------------------------------------------------------------------------------------------------------------------------------------------------------
Using the $k$-Beta transform (\ref{1.17}), we obtain (\ref{2.28})-(\ref{2.29}).
%---------------------------------------------------------------------------------------------------------------------------------------------------------------------------------------------
\end{proof}
%---------------------------------------------------------------------------------------------------------------------------------------------------------------------------------------------
%---------------------------------------------------------------------------------------------------------------------------------------------------------------------------------------------
\begin{theorem}
%---------------------------------------------------------------------------------------------------------------------------------------------------------------------------------------------
The Laplace transform for $\;_{r+1}R_{s,k}$ matrix function are given by
%---------------------------------------------------------------------------------------------------------------------------------------------------------------------------------------------
%---------------------------------------------------------------------------------------------------------------------------------------------------------------------------------------------
\begin{equation}
\begin{split}
&\mathfrak{L}\bigg{[}t^{\frac{C}{k}-I}\;_{r+1}R_{s,k}(A,P_{1},P_{2},\ldots,P_{r};Q_{1},Q_{2},\ldots,Q_{s};B,C,zt^{\frac{B}{k}});\mathbf{s}\bigg{]}\\
&=k(\mathbf{s}k)^{-\frac{C}{k}}\;_{r+1}F_{s,k}(A,P_{1},P_{2},\ldots,P_{r};Q_{1},Q_{2},\ldots,Q_{s};z(k\mathbf{s})^{-\frac{B}{k}}),\label{2.30}
\end{split}
\end{equation}
%---------------------------------------------------------------------------------------------------------------------------------------------------------------------------------------------
\begin{equation}
\begin{split}
&\mathfrak{L}\bigg{[}t^{\frac{E}{k}-I}\;_{r+1}R_{s,k}(A,P_{1},P_{2},\ldots,P_{r};Q_{1},Q_{2},\ldots,Q_{s};B,C,zt);\mathbf{s}\bigg{]}\\
&=k\Gamma_{k}(E)(\mathbf{s}k)^{-\frac{E}{k}}\;_{r+2}R_{s,k}(A,E,P_{1},P_{2},\ldots,P_{r};Q_{1},Q_{2},\ldots,Q_{s};B,C;\frac{z}{\mathbf{s}k}).\label{2.31}
\end{split}
\end{equation}
%---------------------------------------------------------------------------------------------------------------------------------------------------------------------------------------------
\begin{equation}
\begin{split}
&\mathfrak{L}\bigg{[}\;_{r+1}R_{s,k}(A,P_{1},P_{2},\ldots,P_{r};Q_{1},Q_{2},\ldots,Q_{s};B,C;zx);s\bigg{]}\\
&\frac{1}{s}\Gamma_{k}(k)\;_{r+2}R_{s,k}\bigg{(}A,kI,P_{1},P_{2},\ldots,P_{r};Q_{1},Q_{2},\ldots,Q_{s};B,C;\frac{z}{ks}\bigg{)},\label{2.32}
\end{split}
\end{equation}
%---------------------------------------------------------------------------------------------------------------------------------------------------------------------------------------------
\begin{equation}
\begin{split}
&\mathfrak{L}\bigg{[}x^{\frac{A}{k}-I}\;_{r}R_{s,k}(P_{1},P_{2},\ldots,P_{r};Q_{1},Q_{2},\ldots,Q_{s};B,C;zx)\bigg{]}\\
&=\frac{1}{s}(sk)^{I-\frac{A}{k}}\Gamma_{k}(A)\;_{r+1}R_{s,k}\bigg{(}A,P_{1},P_{2},\ldots,P_{r};Q_{1},Q_{2},\ldots,Q_{s};B,C;\frac{z}{ks}\bigg{)}\label{2.33}
\end{split}
\end{equation}
%---------------------------------------------------------------------------------------------------------------------------------------------------------------------------------------------
and
%---------------------------------------------------------------------------------------------------------------------------------------------------------------------------------------------
%---------------------------------------------------------------------------------------------------------------------------------------------------------------------------------------------
\begin{equation}
\begin{split}
&\mathfrak{L}\bigg{[}x^{\frac{P_{i}}{k}-I}\;_{r}R_{s,k}(A,P_{1},P_{2},\ldots,P_{i-1},P_{i+1},\ldots, P_{r};Q_{1},Q_{2},\ldots,Q_{s};B,C;zx)\bigg{]}\\
&=\frac{1}{s}(sx)^{-\frac{P_{i}}{k}}\Gamma_{k}(P_{i})\;_{r+1}R_{s,k}\bigg{(}A,P_{1},P_{2},\ldots,P_{i-1},P_{i},P_{i+1},\ldots,P_{r};Q_{1},Q_{2},\ldots,Q_{s};B,C;\frac{z}{sk}\bigg{)},1\leq i\leq r.\label{2.34}
\end{split}
\end{equation}
%---------------------------------------------------------------------------------------------------------------------------------------------------------------------------------------------
%---------------------------------------------------------------------------------------------------------------------------------------------------------------------------------------------
\end{theorem}
%---------------------------------------------------------------------------------------------------------------------------------------------------------------------------------------------
\begin{proof}
%---------------------------------------------------------------------------------------------------------------------------------------------------------------------------------------------
On applying the Laplace transform (\ref{1.18}), this yields the right hand side of (\ref{2.30})-(\ref{2.34}), we obtain desired results (\ref{2.30})-(\ref{2.34}).
%---------------------------------------------------------------------------------------------------------------------------------------------------------------------------------------------
%---------------------------------------------------------------------------------------------------------------------------------------------------------------------------------------------
\end{proof}
%---------------------------------------------------------------------------------------------------------------------------------------------------------------------------------------------
%---------------------------------------------------------------------------------------------------------------------------------------------------------------------------------------------
\begin{theorem}
%---------------------------------------------------------------------------------------------------------------------------------------------------------------------------------------------
The Fractional $k$-Fourier transform for $\;_{r+1}R_{s,k}$ matrix function is given by
%---------------------------------------------------------------------------------------------------------------------------------------------------------------------------------------------
%---------------------------------------------------------------------------------------------------------------------------------------------------------------------------------------------
\begin{equation}
\begin{split}
&\mathfrak{F}\bigg{[}\;_{r+1}R_{s,k}(A,P_{1},P_{2},\ldots,P_{r};Q_{1},Q_{2},\ldots,Q_{s};B,C,z)\bigg{]}\\
&=\frac{\Gamma_{k}(k)}{iw^{\frac{1}{\alpha}}}\;_{r+2}R_{s,k}\bigg{(}kI,A,P_{1},P_{2},\ldots,P_{r};Q_{1},Q_{2},\ldots,Q_{s};B,C,\frac{i}{kw^{\frac{1}{\alpha}}}\bigg{)}.\label{2.35}
\end{split}
\end{equation}
%---------------------------------------------------------------------------------------------------------------------------------------------------------------------------------------------
%---------------------------------------------------------------------------------------------------------------------------------------------------------------------------------------------
\end{theorem}
%---------------------------------------------------------------------------------------------------------------------------------------------------------------------------------------------
%---------------------------------------------------------------------------------------------------------------------------------------------------------------------------------------------
\begin{proof}
%---------------------------------------------------------------------------------------------------------------------------------------------------------------------------------------------
Applying the Fractional $k$-Fourier transform of $\;_{r+1}R_{s,k}$ for $z<0$, we get
%---------------------------------------------------------------------------------------------------------------------------------------------------------------------------------------------
\begin{equation*}
\begin{split}
&\mathfrak{F}\bigg{[}\;_{r+1}R_{s,k}(A,P_{1},P_{2},\ldots,P_{r};Q_{1},Q_{2},\ldots,Q_{s};B,C,z)\bigg{]}\\
&=\int_{-\infty}^{0}e^{iw^{\frac{1}{\alpha}}z}\;_{r+1}R_{s,k}(A,P_{1},P_{2},\ldots,P_{r};Q_{1},Q_{2},\ldots,Q_{s};B,C,z)dz\\
&=\sum_{\ell=0}^{\infty}\frac{z^{\ell}}{\ell !}(A)_{\ell,k}\prod_{i=1}^{r}(P_{i})_{\ell,k}\bigg{[}\prod_{j=1}^{s}(Q_{j})_{\ell,k}\bigg{]}^{-1}\Gamma^{-1}_{k}(\ell B+C)\int_{-\infty}^{0}e^{iw^{\frac{1}{\alpha}}z}z^{\ell}dz.
\end{split}
\end{equation*}
%---------------------------------------------------------------------------------------------------------------------------------------------------------------------------------------------
%---------------------------------------------------------------------------------------------------------------------------------------------------------------------------------------------
Putting $\chi=-iw^{\frac{1}{\alpha}}z$, $d\chi=-iw^{\frac{1}{\alpha}}dz$ and changing the order of integration and summation, we have
%---------------------------------------------------------------------------------------------------------------------------------------------------------------------------------------------
\begin{equation*}
\begin{split}
&=\sum_{\ell=0}^{\infty}\frac{1}{\ell !}(A)_{\ell,k}\prod_{i=1}^{r}(P_{i})_{\ell,k}\bigg{[}\prod_{j=1}^{s}(Q_{j})_{\ell,k}\bigg{]}^{-1}\Gamma^{-1}_{k}(\ell B+C)\int_{0}^{1}e^{-\chi}\bigg{(}\frac{\chi}{-iw^{\frac{1}{\alpha}}}\bigg{)}^{\ell}\frac{1}{-iw^{\frac{1}{\alpha}}}d\chi\\
&=\sum_{\ell=0}^{\infty}\frac{1}{\ell !}(A)_{\ell,k}\prod_{i=1}^{r}(P_{i})_{\ell,k}\bigg{[}\prod_{j=1}^{s}(Q_{j})_{\ell,k}\bigg{]}^{-1}\Gamma^{-1}_{k}(\ell B+C)(-1)^{\ell}i^{-(\ell+1)}w^{-\frac{\ell+1}{\alpha}}\Gamma(\ell+1)\\
&=\sum_{\ell=0}^{\infty}\frac{1}{\ell !}(A)_{\ell,k}\prod_{i=1}^{r}(P_{i})_{\ell,k}\bigg{[}\prod_{j=1}^{s}(Q_{j})_{\ell,k}\bigg{]}^{-1}\Gamma^{-1}_{k}(\ell B+C)i^{\ell-1}w^{-\frac{\ell+1}{\alpha}}k^{-\ell}\Gamma_{k}(k)(k)_{\ell,k}\\
&=\frac{\Gamma_{k}(k)}{iw^{\frac{1}{\alpha}}}\sum_{\ell=0}^{\infty}\frac{1}{\ell !}(A)_{\ell,k}\prod_{i=1}^{r}(P_{i})_{\ell,k}\bigg{[}\prod_{j=1}^{s}(Q_{j})_{\ell,k}\bigg{]}^{-1}\Gamma^{-1}_{k}(\ell B+C)\bigg{(}\frac{i}{kw^{\frac{1}{\alpha}}}\bigg{)}^{\ell}(k)_{\ell,k}.
\end{split}
\end{equation*}
%---------------------------------------------------------------------------------------------------------------------------------------------------------------------------------------------
\end{proof}
%---------------------------------------------------------------------------------------------------------------------------------------------------------------------------------------------
%---------------------------------------------------------------------------------------------------------------------------------------------------------------------------------------------
%---------------------------------------------------------------------------------------------------------------------------------------------------------------------------------------------
\begin{theorem}
%---------------------------------------------------------------------------------------------------------------------------------------------------------------------------------------------
%---------------------------------------------------------------------------------------------------------------------------------------------------------------------------------------------
Let $\mu>0$, $0<\mu\leq 1$ and $\mathbb{I}_{k}^{\mu}$ be the operators of Riemann-Liouville fractional integral and  $\mathbb{D}_{k}^{\beta}$ fractional derivative then there hold the relations:
%---------------------------------------------------------------------------------------------------------------------------------------------------------------------------------------------
\begin{equation}
\begin{split}
&\mathbb{I}_{k}^{\mu}\bigg{[}t^{\frac{E}{k}}\;_{r+1}R_{s,k}(A,P_{1},P_{2},\ldots,P_{r};Q_{1},Q_{2},\ldots,Q_{s};B,C;t)\bigg{]}(z)\\
&=\Gamma_{k}(E+k I)\Gamma^{-1}_{k}(E+(\mu+k)I)z^{\frac{(\mu+k)I+E}{k}}\\
&\times\;_{r+2}R_{s+1,k}(A,E+kI,P_{1},P_{2},\ldots,P_{r};E+(\alpha+k)I,Q_{1},Q_{2},\ldots,Q_{s};B,C;z),\label{2.36}
\end{split}
\end{equation}
%---------------------------------------------------------------------------------------------------------------------------------------------------------------------------------------------
\begin{equation}
\begin{split}
&\mathfrak{D}_{k}^{\mu}z^{\frac{E}{k}}\;_{r+1}R_{s,k}(A,P_{1},P_{2},\ldots,P_{r};Q_{1},Q_{2},\ldots,Q_{s};B,C,z)\\
&=\frac{1}{k}\Gamma_{k}(E+kI)\Gamma_{k}^{-1}(E+(1-\mu)I)z^{\frac{E+(1-\mu)I}{k}-I}\\
&\times\;_{r+2}R_{s+1,k}(A,E+kI,P_{1},P_{2},\ldots,P_{r};E+(1-\mu)I;B,C,z),Q_{1},Q_{2},\ldots,Q_{s};B,C,z),\label{2.37}
\end{split}
\end{equation}
%---------------------------------------------------------------------------------------------------------------------------------------------------------------------------------------------
%---------------------------------------------------------------------------------------------------------------------------------------------------------------------------------------------
\begin{equation}
\begin{split}
&\mathbb{I}_{a^{+},k}^{\mu}\bigg{[}(t-a)^{\frac{E}{k}-I}\;_{r+1}R_{s,k}(A,P_{1},P_{2},\ldots,P_{r};Q_{1},Q_{2},\ldots,Q_{s};B,C;\nu(t-a))\bigg{]}(z)\\
&=\Gamma_{k}(E)\Gamma^{-1}(E+\mu I)(z-a)^{\frac{E+\mu I}{k}-I}\\
&\times\;_{r+2}R_{s+1,k}(E,A,P_{1},P_{2},\ldots,P_{r};E+\mu I,Q_{1},Q_{2},\ldots,Q_{s};B,C;\nu(z-a))\label{2.38}
\end{split}
\end{equation}
%---------------------------------------------------------------------------------------------------------------------------------------------------------------------------------------------
and
%---------------------------------------------------------------------------------------------------------------------------------------------------------------------------------------------
\begin{equation}
\begin{split}
&\mathbb{D}_{a^{+}}^{\mu}\bigg{[}(t-a)^{\frac{E}{k}-I}\;_{r+1}R_{s,k}(A,P_{1},P_{2},\ldots,P_{r};Q_{1},Q_{2},\ldots,Q_{s};B,C;\nu(t-a))\bigg{]}(z)\\
&=\Gamma_{k}(E)\Gamma_{k}^{-1}(E-\mu I)(z-a)^{\frac{E-\mu I}{k}-I}\\
&\times\;_{r+2}R_{s+1,k}(E,A,P_{1},P_{2},\ldots,P_{r};E-\mu I,Q_{1},Q_{2},\ldots,Q_{s};B,C;\nu(z-a)). \label{2.39}
\end{split}
\end{equation}
%---------------------------------------------------------------------------------------------------------------------------------------------------------------------------------------------
%---------------------------------------------------------------------------------------------------------------------------------------------------------------------------------------------
\end{theorem}
%---------------------------------------------------------------------------------------------------------------------------------------------------------------------------------------------
%---------------------------------------------------------------------------------------------------------------------------------------------------------------------------------------------
\begin{proof}
%---------------------------------------------------------------------------------------------------------------------------------------------------------------------------------------------
To prove assertion (\ref{2.36}). Using equations (\ref{2.1}) and (\ref{1.13}) this yields the right hand side of (\ref{2.36}), we get
%---------------------------------------------------------------------------------------------------------------------------------------------------------------------------------------------
\begin{equation*}
\begin{split}
&\bigg{[}\mathbb{I}_{k}^{\mu}t^{\frac{E}{k}}\;_{r+1}R_{s,k}(A,P_{1},P_{2},\ldots,P_{r};Q_{1},Q_{2},\ldots,Q_{s,k};B,C;t)\bigg{]}(x)\\
&=\frac{1}{k\Gamma_{k}(\mu)}\int_{0}^{z}(z-u)^{\frac{\mu}{k}-1}\;_{r+1}R_{s,k}(A,P_{1},P_{2},\ldots,P_{r};Q_{1},Q_{2},\ldots,Q_{s};B,C;t)t^{\frac{E}{k}}dt\\
&=\sum_{\ell=0}^{\infty}\frac{1}{\ell !}(A)_{\ell,k}\prod_{i=1}^{r}(P_{i})_{\ell,k}\bigg{[}\prod_{j=1}^{s}(Q_{j})_{\ell,k}\bigg{]}^{-1}\Gamma^{-1}_{k}(\ell B+C)\frac{1}{k\Gamma_{k}(\mu)}\int_{0}^{z}(z-t)^{\frac{\mu}{k}-1}t^{\frac{E}{k}+\ell}dt
\end{split}
\end{equation*}
%---------------------------------------------------------------------------------------------------------------------------------------------------------------------------------------------
%---------------------------------------------------------------------------------------------------------------------------------------------------------------------------------------------
Putting $t=z\chi$, $dt=zd\chi$ and changing the order of summation and integration, we obtain (\ref{2.36}). Similarly applying Riemann-Liouville fractional integral and derivative operators, one find the required results (\ref{2.37})-(\ref{2.39}).
%---------------------------------------------------------------------------------------------------------------------------------------------------------------------------------------------
\end{proof}
%---------------------------------------------------------------------------------------------------------------------------------------------------------------------------------------------
%---------------------------------------------------------------------------------------------------------------------------------------------------------------------------------------------
\begin{theorem}
%---------------------------------------------------------------------------------------------------------------------------------------------------------------------------------------------
For the $\;_{r+1}R_{s,k}$ matrix function, we have
%---------------------------------------------------------------------------------------------------------------------------------------------------------------------------------------------
\begin{equation}
\begin{split}
&\mathbb{W}_{k}^{\beta}\bigg{[}(u+a)^{-\frac{E}{k}}\;_{r+1}R_{s,k}\bigg{(}A,P_{1},P_{2},\ldots,P_{r};Q_{1},Q_{2},\ldots,Q_{s};B,C;\frac{1}{u+a}\bigg{)}\bigg{]}(z)\\
&=(z+a)^{\frac{\beta I-E}{k}}\Gamma_{k}^{-1}(E)\Gamma_{k}(E-\beta I)\;_{r+2}R_{s+1,k}\bigg{(}A,E-\beta I,P_{1},P_{2},\ldots,P_{r};E,Q_{1},Q_{2},\ldots,Q_{s};B,C;\frac{1}{z+a}\bigg{)}, \label{2.40}
\end{split}
\end{equation}
%---------------------------------------------------------------------------------------------------------------------------------------------------------------------------------------------
\begin{equation}
\begin{split}
&\mathbb{W}_{k}^{\beta}\bigg{[}(u+a)^{-\frac{C}{k}}\;_{r+1}R_{s,k}\bigg{(}A,P_{1},P_{2},\ldots,P_{r};Q_{1},Q_{2},\ldots,Q_{s};B,C-\beta I;\nu(u+a)^{-\frac{B}{k}}\bigg{)}\bigg{]}\\
&=(z+a)^{\frac{\beta I-C}{k}}\;_{r+1}R_{s,k}\bigg{(}A,P_{1},P_{2},\ldots,P_{r};Q_{1},Q_{2},\ldots,Q_{s};B,C;\nu(z+a)^{-\frac{B}{k}}\bigg{)},  \label{2.41}
\end{split}
\end{equation}
%---------------------------------------------------------------------------------------------------------------------------------------------------------------------------------------------
%---------------------------------------------------------------------------------------------------------------------------------------------------------------------------------------------
\begin{equation}
\begin{split}
&\mathbb{W}_{k}^{-\beta}\bigg{[}\;_{r+1}R_{s,k}\bigg{(}A,P_{1},P_{2},\ldots,P_{r};Q_{1},Q_{2},\ldots,Q_{s,k};B,C;\frac{1}{u+a}\bigg{)}\bigg{]}(z)\\
&=\frac{1}{k}(z+a)^{\frac{(1-\beta)I-E}{k}-I}\Gamma_{k}^{-1}(E)\Gamma_{k}(E+\beta I)\;_{r+2}R_{s+1,k}\bigg{(}A,E+\beta I,P_{1},P_{2},\ldots,P_{r};E,Q_{1},Q_{2},\ldots,Q_{s};B,C;\frac{1}{z+a}\bigg{)} \label{2.42}
\end{split}
\end{equation}
%---------------------------------------------------------------------------------------------------------------------------------------------------------------------------------------------
and
%---------------------------------------------------------------------------------------------------------------------------------------------------------------------------------------------
\begin{equation}
\begin{split}
&\mathbb{W}_{k}^{-\beta}\bigg{[}(u+a)^{-\frac{C}{k}}\;_{r+1}R_{s,k}\bigg{(}A,P_{1},P_{2},\ldots,P_{r};Q_{1},Q_{2},\ldots,Q_{s};B,C+\beta I;\nu(u+a)^{-\frac{B}{k}}\bigg{)}\bigg{]}\\
&=(z+a)^{\frac{(1-\beta)I-C}{k}-I}\;_{r+1}R_{s,k}\bigg{(}A,P_{1},P_{2},\ldots,P_{r};Q_{1},Q_{2},\ldots,Q_{s};B,C;\nu(z+a)^{-\frac{B}{k}}\bigg{)}. \label{2.43}
\end{split}
\end{equation}
%---------------------------------------------------------------------------------------------------------------------------------------------------------------------------------------------
%---------------------------------------------------------------------------------------------------------------------------------------------------------------------------------------------
%---------------------------------------------------------------------------------------------------------------------------------------------------------------------------------------------
\end{theorem}
%---------------------------------------------------------------------------------------------------------------------------------------------------------------------------------------------
%---------------------------------------------------------------------------------------------------------------------------------------------------------------------------------------------
\begin{proof}
%---------------------------------------------------------------------------------------------------------------------------------------------------------------------------------------------
%---------------------------------------------------------------------------------------------------------------------------------------------------------------------------------------------
Applying the $k$-Weyl fractional integral operator (\ref{1.15}), one gets %---------------------------------------------------------------------------------------------------------------------------------------------------------------------------------------------
%---------------------------------------------------------------------------------------------------------------------------------------------------------------------------------------------
%---------------------------------------------------------------------------------------------------------------------------------------------------------------------------------------------
\begin{equation*}
\begin{split}
&\mathbb{W}_{k}^{\beta}\bigg{[}(u+a)^{-\frac{E}{k}}\;_{r+1}R_{s,k}\bigg{(}A,P_{1},P_{2},\ldots,P_{r};Q_{1},Q_{2},\ldots,Q_{s,k};B,C;\frac{1}{u+a}\bigg{)}\bigg{]}(z)\\
&=\frac{1}{k\Gamma_{k}(\beta)}\int_{z}^{\infty}(u-z)^{\frac{\beta}{k}-1}(u+a)^{-\frac{E}{k}}\;_{r+1}R_{s,k}\bigg{(}A,P_{1},P_{2},\ldots,P_{r};Q_{1},Q_{2},\ldots,Q_{s};B,C;\frac{1}{u+a}\bigg{)}du\\
&=\sum_{\ell=0}^{\infty}\frac{z^{\ell}}{\ell !}(A)_{\ell,k}\prod_{i=1}^{r}(P_{i})_{\ell,k}\bigg{[}\prod_{j=1}^{s}(Q_{j})_{\ell,k}\bigg{]}^{-1}\Gamma^{-1}_{k}(\ell B+C)\frac{1}{k\Gamma_{k}(\beta)}\int_{z}^{\infty}(u-z)^{\frac{\beta}{k}-1}(u+a)^{-\ell I-\frac{E}{k}}du.
\end{split}
\end{equation*}
%---------------------------------------------------------------------------------------------------------------------------------------------------------------------------------------------
%---------------------------------------------------------------------------------------------------------------------------------------------------------------------------------------------
Making the change of variables $t=\frac{u-z}{u+a}$, $u=\frac{z+at}{1-t}$, $u=z\Rightarrow t=0$, $u=\infty\Rightarrow t=1$ $du=\frac{z+a}{(1-t)^{2}}dt$, $u+a=\frac{z+a}{1-t}$, $u-z=\frac{t(z+a)}{1-t}$, we get
%---------------------------------------------------------------------------------------------------------------------------------------------------------------------------------------------
%---------------------------------------------------------------------------------------------------------------------------------------------------------------------------------------------
\begin{equation*}
\begin{split}
&(z+a)^{\frac{\beta I-E}{k}}\Gamma_{k}(-\beta)\sum_{\ell=0}^{\infty}\frac{z^{\ell}}{\ell !}(A)_{\ell,k}\prod_{i=1}^{r}(P_{i})_{\ell,k}\bigg{[}\prod_{j=1}^{s}(Q_{j})_{\ell,k}\bigg{]}^{-1}\Gamma^{-1}_{k}(\ell B+C)(z+a)^{-\ell}\Gamma_{k}^{-1}(E)\\
&\times\Gamma_{k}(E-\beta I)(E-\beta I)_{\ell,k}[(E)_{\ell,k}]^{-1}=(z+a)^{\frac{\beta I-E}{k}}\Gamma_{k}^{-1}(E)\Gamma_{k}(E-\beta I)\\
&\times\;_{r+2}R_{s+1,k}\bigg{(}A,E-\beta I,P_{1},P_{2},\ldots,P_{r};E,Q_{1},Q_{2},\ldots,Q_{s,k};B,C;\frac{1}{z+a}\bigg{)}.
\end{split}
\end{equation*}
%---------------------------------------------------------------------------------------------------------------------------------------------------------------------------------------------
%---------------------------------------------------------------------------------------------------------------------------------------------------------------------------------------------
Again applying the $k$-Weyl fractional integral and derivative operators, we get the results (\ref{2.41})-(\ref{2.43}).
%---------------------------------------------------------------------------------------------------------------------------------------------------------------------------------------------
\end{proof}
%---------------------------------------------------------------------------------------------------------------------------------------------------------------------------------------------
%---------------------------------------------------------------------------------------------------------------------------------------------------------------------------------------------
%---------------------------------------------------------------------------------------------------------------------------------------------------------------------------------------------
\section{Some special cases and applications}
%---------------------------------------------------------------------------------------------------------------------------------------------------------------------------------------------
In this section, we develop a integral of the $\;_{r+1}R_{s,k}$ matrix function involving relation with some special cases related to integral representations of $\;_{r+1}R_{s,k}$ matrix functions have also been explained below.
%---------------------------------------------------------------------------------------------------------------------------------------------------------------------------------------------
\begin{theorem} For $|z|<1$, $Re(B)>Re(A)>0$, the $\;_{r+1}R_{r,k}$ matrix function satisfies the following Euler-type
integral representation:
%---------------------------------------------------------------------------------------------------------------------------------------------------------------------------------------------
%---------------------------------------------------------------------------------------------------------------------------------------------------------------------------------------------
\begin{equation}
\begin{split}
\;_{r+1}R_{r,k}(A,\Delta (P,r);\Delta (Q,r);B,C;z)=\Gamma_{k}(Q)\Gamma_{k}^{-1}(P)\Gamma_{k}^{-1}(Q-P)\int_{0}^{1}t^{\frac{P}{k}-I}(1-t)^{\frac{Q-P}{k}-I}\mathbf{E}_{A,B,C}(zt^{\frac{r}{k}})dt\label{3.1}
\end{split}
\end{equation}
%---------------------------------------------------------------------------------------------------------------------------------------------------------------------------------------------
%---------------------------------------------------------------------------------------------------------------------------------------------------------------------------------------------
where $\mathbf{E}_{A,B,C}(z)$ is three parametric Mittag–Leffler matrix function \cite{sd2, snd}.
%---------------------------------------------------------------------------------------------------------------------------------------------------------------------------------------------
\end{theorem}
%---------------------------------------------------------------------------------------------------------------------------------------------------------------------------------------------
%---------------------------------------------------------------------------------------------------------------------------------------------------------------------------------------------
\begin{proof}
%---------------------------------------------------------------------------------------------------------------------------------------------------------------------------------------------
For convenience, let $\;_{r+1}R_{r,k}$ be left-hand side of (\ref{3.1}), then
%---------------------------------------------------------------------------------------------------------------------------------------------------------------------------------------------
\begin{equation}
\begin{split}
&\;_{r+1}R_{r,k}(A,\Delta (P,r);\Delta (Q,r);B,C;z)=\sum_{\ell=0}^{\infty}\frac{z^{\ell}}{\ell !}(A)_{\ell,k}(\frac{1}{r}P)_{\ell,k}(\frac{1}{r}(P+I))_{\ell,k}\ldots\frac{1}{r}(P+(r-1)I)_{\ell,k}\\
&\times[(\frac{1}{r}Q)_{\ell,k}]^{-1}[(\frac{1}{r}(Q+I))_{\ell,k}]^{-1}\ldots[\frac{1}{r}(Q+(r-1)I)_{\ell,k}]^{-1}\Gamma^{-1}_{k}(\ell B+C).\label{3.2}
\end{split}
\end{equation}
%---------------------------------------------------------------------------------------------------------------------------------------------------------------------------------------------
Taking the following properties \cite{sa1}
%---------------------------------------------------------------------------------------------------------------------------------------------------------------------------------------------
%---------------------------------------------------------------------------------------------------------------------------------------------------------------------------------------
\begin{eqnarray}
\begin{split}
(P)_{\ell r,k}=r^{\ell r}\prod_{i=1}^{r}\bigg{(}\frac{P+(i-1)I}{r}\bigg{)}_{\ell,k},\ell=0,1,2,\ldots,\label{3.3}
\end{split}
\end{eqnarray}
%---------------------------------------------------------------------------------------------------------------------------------------------------------------------------------------------
and
%---------------------------------------------------------------------------------------------------------------------------------------------------------------------------------------------
\begin{equation}
\begin{split}
(P)_{r\ell,k}[(Q)_{r\ell,k}]^{-1}=\Gamma_{k}(Q)\Gamma_{k}^{-1}(P)\Gamma_{k}^{-1}(Q-P)\mathbf{B}_{k}(P+rk\ell,Q-P).\label{3.4}
\end{split}
\end{equation}
%---------------------------------------------------------------------------------------------------------------------------------------------------------------------------------------
where $r$ is a positive integer.
%---------------------------------------------------------------------------------------------------------------------------------------------------------------------------------------------
%---------------------------------------------------------------------------------------------------------------------------------------------------------------------------------------------
Using (\ref{3.3}) and (\ref{3.4}) into (\ref{3.2}), we arrive at
%---------------------------------------------------------------------------------------------------------------------------------------------------------------------------------------------
%---------------------------------------------------------------------------------------------------------------------------------------------------------------------------------------------
\begin{equation*}
\begin{split}
&\;_{r+1}R_{r,k}(A,\Delta (P,r);\Delta (Q,r);B,C;z)\\
=&\sum_{\ell=0}^{\infty}\frac{z^{\ell}}{\ell !}(A)_{\ell,k}(\frac{1}{r}P)_{\ell,k}(\frac{1}{r}(P+I))_{\ell,k}(\frac{1}{r}(P+(r-1)I))_{\ell,k}\\
&[(\frac{1}{r}Q)_{\ell,k}]^{-1}[(\frac{1}{r}(Q+I))_{\ell,k}]^{-1}[(\frac{1}{r}(Q+(r-1)I))_{\ell,k}]^{-1}\Gamma^{-1}_{k}(\ell B+C)\\
=&\sum_{\ell=0}^{\infty}\frac{z^{\ell}}{\ell !}(A)_{\ell,k}(P)_{r\ell}[(Q)_{r\ell}]^{-1}\Gamma^{-1}_{k}(\ell B+C)\\
=&\Gamma_{k}(Q)\Gamma_{k}^{-1}(A)\Gamma_{k}^{-1}(Q-P)\sum_{\ell=0}^{\infty}\frac{z^{\ell}}{\ell !}(A)_{\ell,k}\mathbf{B}_{k}(P+r\ell,Q-P)\Gamma^{-1}_{k}(\ell B+C)\\
=&\Gamma_{k}(Q)\Gamma_{k}^{-1}(A)\Gamma_{k}^{-1}(Q-P)\sum_{\ell=0}^{\infty}\frac{z^{\ell}}{\ell !}(A)_{\ell,k}\Gamma^{-1}_{k}(\ell B+C)
\int_{0}^{1}t^{\frac{P+r\ell I}{k}-I}(1-t)^{\frac{Q-P}{k}-I}dt\\
=&\Gamma_{k}(Q)\Gamma_{k}^{-1}(P)\Gamma_{k}^{-1}(Q-P)\int_{0}^{1}t^{\frac{P}{k}-I}(1-t)^{\frac{Q-P}{k}-I}\sum_{\ell=0}^{\infty}\frac{(zt^{\frac{r}{k}})^{\ell}}{\ell !}(A)_{\ell,k}\Gamma^{-1}_{k}(\ell B+C)dt\\
=&\Gamma_{k}(Q)\Gamma_{k}^{-1}(P)\Gamma_{k}^{-1}(Q-P)\int_{0}^{1}t^{\frac{P}{k}-I}(1-t)^{\frac{Q-P}{k}-I}\mathbf{E}_{A,B,C}^{}(zt^{\frac{r}{k}})dt.
\end{split}
\end{equation*}
%---------------------------------------------------------------------------------------------------------------------------------------------------------------------------------------------
\end{proof}
%---------------------------------------------------------------------------------------------------------------------------------------------------------------------------------------------
\begin{theorem}
%---------------------------------------------------------------------------------------------------------------------------------------------------------------------------------------------
For $n,k \in N$, the $\;_{r+1}R_{r,k}$ matrix function satisfy the following Euler type integral representation:
%---------------------------------------------------------------------------------------------------------------------------------------------------------------------------------------------
%---------------------------------------------------------------------------------------------------------------------------------------------------------------------------------------------
\begin{equation}
\begin{split}
&\;_{r+1}R_{r,k}(-n I,\Delta (P,r);\Delta (Q,r);m I,C;z)=\Gamma_{k}(Q)\Gamma_{k}^{-1}(P)\Gamma_{k}^{-1}(Q-P)\\
&\times\Gamma_{k}(n+1)\Gamma_{k}^{-1}(nk I+C)\int_{0}^{1}t^{\frac{P}{k}-I}(1-t)^{\frac{Q}{k}-I}\mathbf{Z}_{n,k}^{C-I}(zt^{\frac{r}{k}};m)dt\label{3.5}
\end{split}
\end{equation}
%---------------------------------------------------------------------------------------------------------------------------------------------------------------------------------------------
%---------------------------------------------------------------------------------------------------------------------------------------------------------------------------------------------
where $\mathbf{Z}_{n,k}^{C-I}(z;m)$ is the Konhauser matrix polynomials \cite{sd1, snd, vct} of degree $n$ in $z^{k}$.
%---------------------------------------------------------------------------------------------------------------------------------------------------------------------------------------------
\end{theorem}
%---------------------------------------------------------------------------------------------------------------------------------------------------------------------------------------------
%---------------------------------------------------------------------------------------------------------------------------------------------------------------------------------------------
\begin{proof}
%---------------------------------------------------------------------------------------------------------------------------------------------------------------------------------------------
By performing $A=-n I$ and $B=m I$ then (\ref{3.1}) reduces to
%---------------------------------------------------------------------------------------------------------------------------------------------------------------------------------------------
%---------------------------------------------------------------------------------------------------------------------------------------------------------------------------------------------
\begin{equation*}
\begin{split}
&\;_{r+1}R_{r,k}(-nI,\Delta (P,r);\Delta (Q,r);m I,C;z)\\
&=\Gamma_{k}(Q)\Gamma_{k}^{-1}(P)\Gamma_{k}^{-1}(Q-P)\int_{0}^{1}t^{\frac{P}{k}-I}(1-t)^{\frac{Q-P}{k}-I}E_{m I,C;-n I}(zt^{r})dt
\end{split}
\end{equation*}
%---------------------------------------------------------------------------------------------------------------------------------------------------------------------------------------------
Using the result defined in \cite{snd}, this leads to right-hand side of (\ref{3.5}).
%---------------------------------------------------------------------------------------------------------------------------------------------------------------------------------------------
\end{proof}
%---------------------------------------------------------------------------------------------------------------------------------------------------------------------------------------------
\begin{theorem}
%---------------------------------------------------------------------------------------------------------------------------------------------------------------------------------------------
For $n\in N$, there reduces to the following integral representation
%---------------------------------------------------------------------------------------------------------------------------------------------------------------------------------------------
%---------------------------------------------------------------------------------------------------------------------------------------------------------------------------------------------
\begin{equation}
\begin{split}
&\;_{r+1}R_{r,k}(-n I,\Delta (P,r);\Delta (Q,r);1 I,C;z)=\Gamma_{k}(Q)\Gamma_{k}^{-1}(P)\Gamma_{k}^{-1}(Q-P)\\
&\times\Gamma_{k}(n+1)\Gamma_{k}^{-1}(Q+nI)\int_{0}^{1}t^{\frac{P}{k}-I}(1-t)^{\frac{Q-P}{k}-I}\mathbf{L}_{n,k}^{C-I}(zt^{\frac{r}{k}})dt,\label{3.6}
\end{split}
\end{equation}
%---------------------------------------------------------------------------------------------------------------------------------------------------------------------------------------------
where $\mathbf{L}_{n,k}^{C-I}(z)$ is the Laguerre matrix polynomials \cite{j1}.
%---------------------------------------------------------------------------------------------------------------------------------------------------------------------------------------------
\end{theorem}
%---------------------------------------------------------------------------------------------------------------------------------------------------------------------------------------------
\subsection{Special cases}
%-----------------------------------------------------------------------------------------------------------------------------------------------------------------
%---------------------------------------------------------------------------------------------------------------------------------------------------------------------------------------
Case 1. On setting $k=1$, $P_{i}=I$, and $Q_{j}=I$, (\ref{2.1}) becomes
%-----------------------------------------------------------------------------------------------------------------------------------------------------------------
\begin{equation}
\begin{split}
\;_{1}R_{0}(A;-;B,C;z)=\sum_{\ell=0}^{\infty}\frac{z^{\ell}}{\ell !}(A)_{\ell}\Gamma^{-1}(\ell B+C)=E_{B,C}^{A}(z),\label{3.7}
\end{split}
\end{equation}
%---------------------------------------------------------------------------------------------------------------------------------------------------------------------------------------
where $E_{B,C}^{A}(z)$ is the three parameter Mittag--Leffler matrix function given by Sanjhira et al. \cite{snd}

Case 2. On setting $k=1$, (\ref{2.1}) reduces to
%-----------------------------------------------------------------------------------------------------------------------------------------------------------------
\begin{equation}
\begin{split}
&\;_{r+1}R_{s,1}(A,P_{1},P_{2},\ldots,P_{r};Q_{1},Q_{2},\ldots,Q_{s};B,C;z)\\
=&\;_{r}K_{s}(A,P_{1},P_{2},\ldots,P_{r};Q_{1},Q_{2},\ldots,Q_{s};B,C;z),\label{3.8}
\end{split}
\end{equation}
%---------------------------------------------------------------------------------------------------------------------------------------------------------------------------------------
where $\;_{r}K_{s}$ is the matrix $K$-function given by Sharma and Jain \cite{sj}.

Case 3. By taking $k=1$ and $A=I$, (\ref{2.1}) becomes
%-----------------------------------------------------------------------------------------------------------------------------------------------------------------
\begin{equation}
\begin{split}
&\;_{r+1}R_{s,1}(I,P_{1},P_{2},\ldots,P_{r};Q_{1},Q_{2},\ldots,Q_{s};B,C;z)\\
&=\;_{r}M_{s}(P_{1},P_{2},\ldots,P_{r};Q_{1},Q_{2},\ldots,Q_{s};B,C;z),\label{3.9}
\end{split}
\end{equation}
%---------------------------------------------------------------------------------------------------------------------------------------------------------------------------------------
where $\;_{r}M_{s}$ is the generalized matrix $M$-series defined by Sharma and Jain \cite{sj}.

Case 4. On setting $k=1$, $A=I$ and $Q_{1}=I$, (\ref{2.1}) reduces to
%-----------------------------------------------------------------------------------------------------------------------------------------------------------------
\begin{equation}
\begin{split}
&\;_{r+1}R_{s,1}(I,P_{1},P_{2},\ldots,P_{r};I,Q_{2},\ldots,Q_{s};B,C;z)\\
=&\;_{r}R_{s-1}(P_{1},P_{2},\ldots,P_{r};Q_{2},\ldots,Q_{s};B,C;z),\label{3.10}
\end{split}
\end{equation}
%---------------------------------------------------------------------------------------------------------------------------------------------------------------------------------------
where $\;_{r}R_{s-1}$ matrix function defined by Sanjhira and Dwivedi \cite{sd2} and Shehata et al. \cite{skc}.

Case 5. On taking $k=1$, $A=I$ and $B=C=I$, (\ref{2.1}) becomes
%-----------------------------------------------------------------------------------------------------------------------------------------------------------------
\begin{equation}
\begin{split}
&\;_{r+1}R_{s,1}(I,P_{1},P_{2},\ldots,P_{r};Q_{1},Q_{2},\ldots,Q_{s};I,I;z)\\
=&\;_{r}F_{s}(P_{1},P_{2},\ldots,P_{r};Q_{1},Q_{2},\ldots,Q_{s};z),\label{3.11}
\end{split}
\end{equation}
%---------------------------------------------------------------------------------------------------------------------------------------------------------------------------------------
where $\;_{r}F_{s}(z)$ is the generalized hypergeometric matrix function by Shehata \cite{sa1}.
%---------------------------------------------------------------------------------------------------------------------------------------------------------------------------------------------
%---------------------------------------------------------------------------------------------------------------------------------------------------------------------------------------------
\section{Conclusion}
%---------------------------------------------------------------------------------------------------------------------------------------------------------------------------------------------
It should be observed in conclusion, we are motivated to obtain the our interest recurrence matrix relation, differential properties, new integral representations, Fractional $k$-Beta, Laplace and $k$-Fourier transforms, the $k$-Riemann–Liouville and $k$-Weyl fractional integral and derivative operators of this $\;_{r+1}R_{s,k}$ matrix function, which is an important crucial in theory of classical analysis, integral transforms, mathematical analysis, operational techniques, mathematical physics, fractional calculus, applied mathematics and statistics.  Hence, for our purposes, the results appear in this study are seemed novel to the literature. From this view, we gives some special cases, such as hypergeometric and Mittag-Leffler matrix functions, $k$-Konhauser and $k$-Laguerre matrix polynomials for the $\;_{r+1}R_{s,k}$ matrix function using transform method with its application to the Mittag–Leffler matrix function and Euler-type integral representations including many special cases and so on. Lastly, we conclude this paper by hoping that we will be able also extend these generalized type $k$-fractional derivatives and their consequences for by analytical continuation. We have been asked to publish in the present form in the hope that others, more qualified than ourselves in this field, may find the results suggestive, even indicative, of the usefulness of a more systematic study.
%-------------------------------------------------------------------------------------------------------------------------------------------------------------------------------------
\begin{definition}
%-------------------------------------------------------------------------------------------------------------------------------------------------------------------------------------
Our definitions of extended $\;_{r+1}R_{s,k}$ matrix function can further results be generalized Pochhammer symbol $(A,\mu)_{\ell,k}$ to the following form:
%---------------------------------------------------------------------------------------------------------------------------------------------------------------------------------------------
\begin{equation}
\begin{split}
&\;_{r+1}R_{s,k}((A,\mu),P_{1},P_{2},\ldots,P_{r};Q_{1},Q_{2},\ldots,Q_{s};B,C;z)\\
=&\sum_{\ell=0}^{\infty}\frac{z^{\ell}}{\ell !}(A,\mu)_{\ell,k}\prod_{i=1}^{r}(P_{i})_{\ell,k}\bigg{[}\prod_{j=1}^{s}(Q_{j})_{\ell,k}\bigg{]}^{-1}\Gamma^{-1}_{k}(\ell B+C),
\end{split}
\end{equation}
%---------------------------------------------------------------------------------------------------------------------------------------------------------------------------------------------
where
%---------------------------------------------------------------------------------------------------------------------------------------------------------------------------------------------
\begin{equation}
\begin{split}
(A,\mu)_{\ell,k}=\Gamma_{k}^{-1}(A)\int_{0}^{\infty}t^{\frac{A+\ell}{k}-1}e^{-t-\frac{\mu}{t}}dt;Re(A)>0, \mu>0
\end{split}
\end{equation}
%---------------------------------------------------------------------------------------------------------------------------------------------------------------------------------------------
%---------------------------------------------------------------------------------------------------------------------------------------------------------------------------------------------
and the extension of extended $\;_{r+1}R_{s,k}$ matrix function will be defined as follows
%---------------------------------------------------------------------------------------------------------------------------------------------------------------------------------------------
\begin{equation}
\begin{split}
&\;_{r+1}R_{s,k}((A,\mu),P_{1},P_{2},\ldots,P_{r};Q_{1},Q_{2},\ldots,Q_{s};B,C;\tau;z)\\
=&\prod_{i=1}^{r}\prod_{j=1}^{s}\Gamma_{k}^{-1}(P_{i})\Gamma_{k}(Q_{j})\sum_{\ell=0}^{\infty}\frac{z^{\ell}}{\ell !}(A,\mu)_{\ell,k}\Gamma_{k}(P_{i}+k\ell\tau)\Gamma_{k}^{-1}(Q_{j}+k\ell\tau)\Gamma^{-1}_{k}(\ell B+C).
\end{split}
\end{equation}
%---------------------------------------------------------------------------------------------------------------------------------------------------------------------------------------------
%-------------------------------------------------------------------------------------------------------------------------------------------------------------------------------------
\end{definition}
%---------------------------------------------------------------------------------------------------------------------------------------------------------------------------------------------
\section*{Funding}This research received no external funding.
%---------------------------------------------------------------------------------------------------------------------------------------------------------------------------------------------
\section*{Data availability}Data sharing is not applicable to this article as no data sets were generated or analyzed during the current study.
%---------------------------------------------------------------------------------------------------------------------------------------------------------------------------------------------
\section*{Declarations}
%---------------------------------------------------------------------------------------------------------------------------------------------------------------------------------------------
\subsection*{Conflicts of Interest}The author declare no conflict of interest.
%---------------------------------------------------------------------------------------------------------------------------------------------------------------------------------------------
%---------------------------------------------------------------------------------------------------------------------------------------------------------------------------------------------
\section*{Ethical approval} This article does not contain any studies with human participants performed by any of the
author.
%---------------------------------------------------------------------------------------------------------------------------------------------------------------------------------------------

%---------------------------------------------------------------------------------------------------------------------------------------------------------------------------------------------
\end{document}